\documentclass[11pt,twoside]{article}
\usepackage[english]{babel}
\usepackage[utf8]{inputenc}
\usepackage{amsmath}
\usepackage{amssymb}
\usepackage{amsfonts}
\usepackage{graphicx}

\usepackage[textwidth=16.2cm,textheight=22cm,top=3cm,bottom=3cm]{geometry}
\usepackage{fancyhdr}

\usepackage{color}
\graphicspath{{./figuras/}}% directorio de figuras
\newtheorem{theorem}{Theorem}
\newtheorem{lemma}{Lemma}

\newcommand{\drop}[1]{}
\newcommand{\no}{\noindent}

\newcommand{\p}{{\partial}}
\newcommand{\eps}{\varepsilon}
\newcommand{\vfi}{\varphi}
\newcommand{\R}{\mathbb{R}}

\newcommand{\fer}[1]{(\ref{#1})}
\newcommand{\qtext}[1]{\quad\text{#1}}
\newcommand{\qtextq}[1]{\quad\text{#1}\quad}
\newcommand{\abs}[1]{| #1 |}
\newcommand{\nor}[1]{\| #1 \|}
\newcommand{\grad}{\nabla}

\newcommand{\bu}{\mathbf{u}}
\newcommand{\bv}{\mathbf{v}}

\newcommand{\bd}{\mathbf{d}}

\newcommand{\bF}{\mathbf{F}}
\newcommand{\bR}{\mathbf{R}}

\newcommand{\bw}{\mathbf{w}}
\newcommand{\bx}{\mathbf{x}}
\newcommand{\by}{\mathbf{y}}
\newcommand{\bG}{\mathbf{G}}

\newcommand{\cN}{\mathcal{N}}
\newcommand{\cL}{\mathcal{L}}
\newcommand{\cQ}{\mathcal{Q}}
\newcommand{\cM}{\mathcal{M}}
\newcommand{\cS}{\mathcal{S}}
\newcommand{\cA}{\mathcal{A}}

\newcommand{\cR}{\mathcal{R}}

\renewcommand{\bf}{\mathbf{f}}

\def\O{\Omega}
\def\brho{\mbox{\boldmath$\rho$}}
\def\bpsi{\mbox{\boldmath$\psi$}}
\def\boldeta{\mbox{\boldmath$\eta$}}

\DeclareMathOperator{\Div}{div}
\DeclareMathOperator{\tr}{tr}

\DeclareMathOperator{\RD}{RD}

\title{Turing instability analysis of a singular cross-diffusion problem
\thanks{First author supported by the Spanish MEC Project MTM2017-87162-P.
}}
\author{Gonzalo Galiano\thanks{Dept. of Mathematics, University of Oviedo, Spain 
  ({\tt galiano@uniovi.es}).}
\and V\'{\i}ctor Gonz\'alez-\!Tabernero\footnotemark[3] \thanks{Universidad de Santiago de Compostela, Santiago de Compostela, Spain ({\tt victor.gonzalez.tabernero@rai.usc.es}).}}
  
\date{}

\begin{document}

\maketitle

\begin{abstract}
%Turing instability is a well known mechanism of pattern formation in diffusive systems.
The population model of Busenberg and Travis is a paradigmatic model in ecology and tumour modelling due to its ability to capture interesting phenomena like the segregation of populations. Its singular mathematical structure enforces the consideration of regularized problems to deduce properties as fundamental as the existence of solutions. In this article we perform a weakly nonlinear stability analisys of a general class of regularized problems to study the convergence of the instability modes in the limit of the regularization parameter. We demonstrate with some specific examples that the pattern formation observed in the regularized problems, with unbounded wave numbers, is not present in the limit problem due to the amplitude decay of the oscillations. We also check the results of the stability analysis with direct finite element simulations of the problem. 

\end{abstract}

\no\emph{Keywords: }
Cross-diffusion, Turing instability, weakly nonlinear analysis, finite element.

In \cite{Busenberg83}, Busenberg and Travis introduced a class of singular cross-diffusion problems  under the assumption that the spatial relocation of each species is due to a diffusion flow which depends on the densities of all the involved species. In the case of two species, if $u_1$, $u_2$ denote their densities, the flow, in its simplest form, may be assumed to be determined by the total population $u_1+u_2$, and thus the conservation laws for both species lead to the  system	
\begin{align}
&  \p_t u_1 - \Div \big( u_1 (\grad u_1 +\grad u_2)\big) = f_1(u_1,u_2),  \label{eq:u1:BT} \\
&  \p_t u_2 -  \Div \big( u_2 (\grad u_1 +\grad u_2)\big) =  f_2(u_1,u_2). \label{eq:u2:BT}
\end{align}
The functions $f_1$ and $f_2$ capture some ecological features of the populations, such as growth, competition, etc. As usual, the equations \fer{eq:u1:BT}-\fer{eq:u2:BT} are complemented with non-negative initial data and non-flow boundary conditions.

The system \fer{eq:u1:BT}-\fer{eq:u2:BT} is called a \emph{cross-diffusion} system because the flow of each species depend upon the densities of the other species. We call it \emph{singular} because the resulting diffusion matrix is singular. Indeed, when rewritting \fer{eq:u1:BT}-\fer{eq:u2:BT} in matrix form, for $\bu=(u_1,u_2)$, we get the equation
\begin{align*}
%\label{def:diffmat}
\p_t\bu -\Div(\cA(\bu)\grad\bu)=\bf(\bu), \qtextq{with} \cA(\bu)=\begin{pmatrix} u_1& u_1\\ u_2&  u_2 \end{pmatrix},	
\end{align*}
where the divergence is applied by rows. 
The full and singular structure of $\cA$ introduces serious difficulties in the mathematical analysis of the problem, as we shall comment later.

In his seminal paper \cite{Turing1952}, Turing introduced a mechanism explaining how spatially uniform equilibria may evolve, small perturbations mediating, into stable equilibria with non-trivial spatial  structure. He considered a system of the type 
\begin{align}
&  \p_t u_1 - \Delta u_1 = f_1(u_1,u_2),  \label{eq:u1:ld}\\
&  \p_t u_2 - \sigma \Delta u_2 =  f_2(u_1,u_2),\label{eq:u2:ld}
\end{align}
with $\sigma>0$, and proved that when $\sigma$ is small or large enough then the stable equilibria of the dynamical system 
\begin{align}
&  \p_t v_1  = f_1(v_1,v_2),  \label{eq:u1:DS}\\
&  \p_t v_2  =  f_2(v_1,v_2),\label{eq:u2:DS}
\end{align}
are not stable for the diffusion system \fer{eq:u1:ld}-\fer{eq:u2:ld} and that, in their place,  non-uniform equilibria with spatial structure become the new stable solutions. This mechanism is known as \emph{Turing instability} or \emph{Turing bifurcation}.

In this article we study Turing instability for the cross-diffusion singular system \fer{eq:u1:BT}-\fer{eq:u2:BT}. We already know that some cross-diffusion systems, such as the paradigmatic SKT model introduced by Shigesada, Kawasaki and Teramoto \cite{Shigesada1979}, exhibit Turing instability when cross-diffusion coefficients are large in comparison with self-diffusion coefficients, see e.g. \cite{Gambino2012, Gambino2013}. However, the singularity of the diffusion matrix of the system \fer{eq:u1:BT}-\fer{eq:u2:BT}, not present in the SKT model,  introduces important mathematical difficulties to the analysis of this system. 

Regarding the existence of solutions of \fer{eq:u1:BT}-\fer{eq:u2:BT}, it has been proved only in some special situations:  for a bounded spatial domain $\O\subset\R$  (Bertsch et al. \cite{Bertsch2010}) and for  $\O= \R^n$ (Bertsch et al. \cite{Bertsch2012}). In their proofs, the following observation is crucial: adding the two equations of \fer{eq:u1:BT}-\fer{eq:u2:BT} shows that if a solution of this system does exist then the total population, $u=u_1+u_2$, satisfies the porous medium type equation 
\begin{align}
\p_t u -\Div(u\grad u) = f(u), 	\label{eq:PM}
\end{align}
for which the theory of existence and uniqueness of solutions is well established. In particular, if the initial data of the total population is bounded away from zero and if $f$ is regular enough with $f(0)\geq0$, it is known that the solution of \fer{eq:PM} remains positive and smooth for all time. This allows to introduce the change of unknowns $w_i=u_i/u$, for $i=1,2$, into the original problem \fer{eq:u1:BT}-\fer{eq:u2:BT} to deduce the equivalent formulation 
\begin{align}
& \p_t u -\Div(u\grad u) = F_1(u,w_1), \label{eq:u:B}\\	
& \p_t w_1 - \grad u \cdot \grad w_1 = F_2(u,w_1),\label{eq:w1:B}
\end{align}
for certain well-behaved functions $F_1$ and $F_2$. Being the structure of the system \fer{eq:u:B}-\fer{eq:w1:B} of parabolic-hyperbolic nature, parabolic regularization of the system by adding the term $-\delta \Delta w_1$ to the left hand side of \fer{eq:w1:B}, and the consideration of the characteristics defined by the field $\grad u$ are the main ingredients of the proofs made by Bertsch et al. \cite{Bertsch2010, Bertsch2012}. 

In \cite{Galiano2015} we followed a different approach to prove the existence of solutions of the original system \fer{eq:u1:BT}-\fer{eq:u2:BT} for a bounded domain $\O\subset\R$. We directly performed a parabolic regularization of the system by introducing a cross-diffusion perturbation term while keeping the porous medium type equation satisfied by $u$. More concretely, we considered the system 
\begin{align}
&  \p_t u_1 - \Div \big( u_1 (\grad u_1 +\grad u_2)\big)-\frac{\delta}{2}\Delta(u_1(u_1+u_2)) = f_1(u_1,u_2),  \label{eq:u1:GS} \\
&  \p_t u_2 - \Div \big( u_2 (\grad u_1 +\grad u_2)\big) -\frac{\delta}{2}\Delta(u_2(u_1+u_2)) =  f_2(u_1,u_2), \label{eq:u2:GS}
\end{align}
and then used previous results for cross-diffusion systems \cite{Galiano2003, Chen2004, Jungel2015} to establish the existence of solutions of the approximated problems. Then, BV estimates  similar to those obtained in \cite{Bertsch2010} allowed to prove the convergence of the sequence $(u_1^{(\delta)},u_2^{(\delta)} )$ to a solution of the original problem.  Let us finally mention that the system \fer{eq:u1:BT}-\fer{eq:u2:BT} is a limit case of a general type of problems with diffusion matrix given by
\begin{align*}
	\cA(\bu)=\begin{pmatrix} a_{11}u_1& a_{12}u_1\\ a_{21}u_2& a_{22} u_2 \end{pmatrix},	
\end{align*}
for which, if $a_{ii}>0$, for $i=1,2$, and $a_{11}a_{22}>a_{12}a_{21}$ then the existence of solutions in ensured for any spatial dimension of $\O$, see \cite{Galiano2014}. In addition, it has been shown that this kind of systems, when set in the whole space $\O=\R^n$, may be obtained as mean field limits \cite{Chen2019}.  

 Concerning Turing instability, since the diffusion matrix, $\cA(\bu)$, corresponding to the system \fer{eq:u1:BT}-\fer{eq:u2:BT} is singular, the linearization of this system about an equilibrium of the dynamical system \fer{eq:u1:DS}-\fer{eq:u2:DS} does not provide any information on the behaviour of the equilibrium in the spatial dependent case.  Thus, our approach to the investigation of Turing instability for the system \fer{eq:u1:BT}-\fer{eq:u2:BT} relies on the study of this property for approximating problems like \fer{eq:u1:GS}-\fer{eq:u2:GS} and its limit behaviour.

% Fixing $\delta>0$, it is not difficult to prove that instability arises for extreme values of $\sigma$. This is not surprising after the results found for simpler systems like \fer{eq:u1:ld}-\fer{eq:u2:ld}. Thus, for concreteness, we focus in the case $\sigma=1$ and show the occurrence of Turing instability when $\delta$ is small enough, which corresponds to the more interesting case of cross-diffusion driven instability.  
 
We prove that linear instability is always present in the limit $\delta\to0$, which is the case  when the sequence of solutions of the approximated problems \fer{eq:u1:GS}-\fer{eq:u2:GS} converges to the solution of the original problem  \fer{eq:u1:BT}-\fer{eq:u2:BT}. Interestingly, the linear analysis also establishes that the main instability wave number is unbounded as $\delta\to0$.
 
 For a clearer understanding of this convergence of a increasingly oscillating sequence of functions to a $BV$ function (the solution of  \fer{eq:u1:BT}-\fer{eq:u2:BT} ensured in \cite{Bertsch2010, Galiano2015}), we perform a weakly nonlinear analysis (WNA) which allows to gain insight into the behaviour of the amplitude of the main instability mode as  $\delta\to0$. As expected, we find that the amplitude of the instability modes vanishes in the limit $\delta\to0$ resulting, therefore, coherent with the $BV$ convergence. In addition, this result also suggests that the uniform equilibrium is stable for the original problem. We furthermore check these analytical results by numerically comparing the WNA approximation to a FEM approximation of the nonlinear problem.

\section{Main results}

For simplicity, we study Turing instability for the one-dimensional spatial setting which has also the advantage of a well stablished existence theory for the case of a bounded domain \cite{Bertsch2010, Galiano2015}. By redefining the functions $f_1,f_2$, we can fix without loss of generality $\O=(0,\pi)$ and then rewrite problem \fer{eq:u1:BT}-\fer{eq:u2:BT} together with the usual auxiliary conditions as 
\begin{align}
&  \p_t u_1 - \p_x \big( u_1 (\p_x u_1 +\p_x u_2)\big) =  f_1(u_1,u_2) & &\text{in } Q_T, & \label{eq:u1}\\
&  \p_t u_2 -  \p_x \big( u_2 (\p_x u_1 +\p_x u_2)\big) = f_2(u_1,u_2) & &\text{in } Q_T ,\label{eq:u2}&\\
& u_1 (\p_x u_1 +\p_x u_2) = u_2 (\p_x u_1 +\p_x u_2)   = 0& &\text{on }  \Gamma_T, & \label{eq:bd}\\
& u_1(0,\cdot) = u_{10}, \quad u_2(0,\cdot) = u_{20} & &\text{in } \O, & \label{eq:id}
\end{align}
where $Q_T=(0,T)\times\O$ and the initial data $u_{10},u_{20}$ are non-negative functions. We assume a competitive Lotka-Volterra form for the reaction term, this is,
$f_i(u_1,u_2) = u_i(\alpha_i-\beta_{i1}u_1-\beta_{i2}u_2)$,	
 for $i=1,2$, and for some non-negative parameters $\alpha_i,\beta_{ij}$, for $i,j=1,2$.

In order to deal with several types of regularized problems we introduce, for positive $\delta$ and $b$, the  uniformly parabolic cross-diffusion system
\begin{align}
&  \p_t u_1 - \p_x \big( d_{11}^\delta(\bu)\p_x u_1 +d_{12}^\delta(\bu)\p_x u_2\big) =  f_1^b(\bu) & &\text{in } Q_T, & \label{eq:u1g}\\
&  \p_t u_2 - \p_x \big( d_{21}^\delta(\bu)\p_x u_1 +d_{22}^\delta(\bu)\p_x u_2\big) =  f_2^b(\bu) & &\text{in } Q_T ,\label{eq:u2g}&\\
& d_{11}^\delta(\bu)\p_x u_1 +d_{12}^\delta(\bu)\p_x u_2 = d_{21}^\delta(\bu)\p_x u_1 +d_{22}^\delta(\bu)\p_x u_2   = 0& &\text{on }  \Gamma_T, & \label{eq:bdg}\\
& u_1(0,\cdot) = u_{10}, \quad u_2(0,\cdot) = u_{20} & &\text{in } \O, & \label{eq:idg}
\end{align}
where the diffusion matrix $D^\delta(\bu) = (d_{ij}^\delta(\bu))$ and the Lotka-Volterra function $\bf^b(\bu) =(f_1^b(\bu),f_2^b(\bu))$  satisfy the assumptions $H_D$:
\begin{enumerate}
\item $D^\delta(\bu)$ is linear in $\bu$ and affine in $\delta$, so that it allows the decompositions 
\begin{align}
\label{decomp:Ddelta}
D^\delta(\bu) = D^{0}(\bu)+\delta D^{1}(\bu)	 = D^{\delta 1}u_1 + D^{\delta 2}u_2 ,
\end{align}
for some matrices $D^{\delta i}$ for $i=1,2$, being the coefficients of $D^\delta(\bu)$ given by 
\begin{align*}
%\label{def:dijdelta}
 d_{ij}^\delta(\bu) = d_{ij}^{10}u_1+ d_{ij}^{11}u_1 \delta + d_{ij}^{20}u_2 + d_{ij}^{21}u_2 \delta , %=\sum_{m=1}^2\sum_{n=0}^1 d_{ij}^{mn}u_m\delta^n,
\end{align*}
for some non-negative constants $d_{ij}^{m n}$, for $i,j,m=1,2$ and $n=0,1$.  
\item We assume that $d_{ii}^\delta(\bu)>0$ for $i=1,2$, and that $\det(D^\delta(\bu))$ is an increasing function with respect to $\delta$ satisfying $\det(D^\delta(\bu))>0$ if $\delta>0$ and $\bu\in\R_+^2$.
\item For $i,j=1,2$, $f_i^b(u_1,u_2) = u_i(\alpha_i^b-\beta_{i1}^bu_1-\beta_{i2}^bu_2)$ for some non-negative $\alpha_i^b,\beta_{ij}^b$  such that $\alpha_i^b\to\alpha_i$ and $\beta_{ij}^b\to\beta_{ij}$ as $b\to0$. Moreover, using the notation $\alpha_i^0=\alpha_i$ and $\beta_{ij}^0=\beta_{ij}$, we assume, for $b\geq0$,
\begin{align}
\label{ass:B}
\begin{split}
    &\beta_{22}^b\alpha_1^b - \beta_{12}^b\alpha_2^b > 0, \hspace{3mm} \beta_{11}^b\alpha_2^b - \beta_{21}^b\alpha_1^b > 0, \\
    &\det(B^b)>0, \hspace{3mm} \tr(B^b)\geq 0, \hspace{3mm} \text{where } B^b=\left( \beta_{ij}^b\right).
\end{split}
\end{align}

\end{enumerate}
Observe that \fer{ass:B} ensures the existence of a stable coexistence equilibrium for the dynamical system  \fer{eq:u1:DS}-\fer{eq:u2:DS}, given by   
\begin{equation}
\label{equilibrium}
    \mathbf{u}^*=\left(u_1^*, u_2^*\right) = \left( \frac{\beta_{22}^b\alpha_1^b - \beta_{12}^b\alpha_2^b}{\beta_{11}^b \beta_{22}^b- \beta_{12}^b\beta_{21}^b},\frac{\beta_{11}^b\alpha_2^b - \beta_{21}^b\alpha_1^b}{\beta_{11}^b \beta_{22}^b- \beta_{12}^b\beta_{21}^b} \right).
\end{equation}

%Observe that since $D^\delta(\bu)$ is defined to be a smooth approximation of the diffusion matrix $\cA(\bu)$ of the system \fer{eq:u1}-\fer{eq:u2}, see \fer{def:diffmat}, we necessarily must define 
%\begin{align*}
%D^0(\bu)=\begin{pmatrix}
%u_1 & u_1\\ u_2 & u_2	
%\end{pmatrix}.
%\end{align*}
There are two examples of $D^\delta(\bu)$ in which we are specially interested. The first, due to its simplicity for the calculations. 
We set
\begin{align}
\label{def:ex1}
D^\delta(\bu) =
\begin{pmatrix}
 (1+\delta)u_1 & u_1\\
  u_2 & (1+\delta)u_2
 \end{pmatrix},
\end{align}
for which $\det(D^\delta(\bu)) = \delta(2+\delta)u_1u_2$.
According to \cite{Galiano2014}, the second hypothesis of $H_D$ guarantees the well-posedness of the problem \fer{eq:u1g}-\fer{eq:idg} corresponding to this diffusion matrix. 
The second example corresponds to the approximation used in \cite{Galiano2015} for proving the existence of $BV$ solutions  of the original problem \fer{eq:u1}-\fer{eq:id}:
%.  We consider the diffusion matrix of the system \fer{eq:u1:GS}-\fer{eq:u2:GS}:
\begin{align}
\label{def:ex2}
D^\delta(\bu) =
\begin{pmatrix}
 (1+\delta)u_1+\frac{\delta}{2}u_2 & (1+\frac{\delta}{2})u_1\\[0.25em]
  (1+\frac{\delta}{2})u_2 & \frac{\delta}{2}u_1 + (1+\delta)u_2
 \end{pmatrix},
\end{align}
for which $\det(D^\delta(\bu)) = \frac{1}{2}\delta(1+\delta)(u_1+u_2)^2$. 

The approximation of the reaction terms introduced in the system \fer{eq:u1g}-\fer{eq:idg} is not essential. Its aim is to support the specific example we deal with in Theorem~\ref{th:example}, but can be ignored ($b=0$) in the general linear and weakly nonlinear analysis of Theorems~\ref{th:estab} and \ref{th:wna}.  Nevertheless, we state these results taking it into account.
Our first result gives conditions under which linear instability arises. 
The following notation is used:
\begin{align}
\label{def:K}
    K=D\bf^b(\bu^*)=\begin{pmatrix}
    -\beta_{11}^b u_1^* & -\beta_{12}^b u_1^*\\
    -\beta_{21}^b u_2^* & -\beta_{22}^b u_2^*\\
    \end{pmatrix}.
\end{align}
\begin{theorem}[Linear instability]
\label{th:estab}
Assume $H_D$, with $b\geq0$. Let $\bu^*$ be the coexistence equilibrium defined by  \fer{equilibrium}. If%, for all $\delta\geq0 $, we have  
\begin{align}
\label{ass:inst2}
\tr(K^{-1} D^\delta(\bu^*))>0 	\qtextq{for all}\delta\geq0
\end{align}
then there exists $\delta_c>0$ such that if $\delta<\delta_c$ then 
 $\bu^*$ is a linearly unstable equilibrium for problem  \fer{eq:u1g}-\fer{eq:idg}. 
  In such situation, the wave number of the main instability mode tends to infinity as $\delta\to0$.
\end{theorem}
Condition \fer{ass:inst2} is equivalent to 
\begin{align}
\label{ass:inst}
 d_{11}^\delta(\bu^*)\beta_{22}^bu_2^*+ d_{22}^\delta(\bu^*)\beta_{11}^bu_1^* < d_{12}^\delta(\bu^*)\beta_{21}^bu_2^*+ d_{21}^\delta(\bu^*)\beta_{12}^bu_1^*
\end{align}
and introduces a further restriction on the matrix of competence coefficients. Roughly speaking, for $B^b$ to fulfil both \fer{ass:B} and \fer{ass:inst}, its elements must be such that intra-population joint competence is larger than inter-population joint competence (condition \fer{ass:B}) and one of the inter-population competence coefficients is large in comparison with the others (condition \fer{ass:inst}). A numeric example we shall work with along the article is
\begin{align}
\label{def:Bex}
B^b=\begin{pmatrix}
1& \dfrac{b}{2}\\ 2 & 1	
\end{pmatrix}, \qtextq{with } b\in(0,\tfrac{1}{2}).
\end{align}
Assuming the forms of $D^\delta(\bu^*)$ given in Examples 1 and 2, see \fer{def:ex1} and \fer{def:ex2}, we have that  the conditions   \fer{ass:B} and \fer{ass:inst} on $B^b$ are satisfied if  $\delta<b/4$  (Example 1) or $\delta<bu_1^*u_2^*/(u_1^*+u_2^*)^2$ (Example 2). Therefore, the most meaningful case when  $\delta$ is close to zero is satisfied by both diffusion matrices.

\bigskip

%The next result gives sufficient conditions to have the interval $(0,1]$ of admissible values of $\sigma$ split in a region of instability, $(0,\sigma_c)$,  and a region of stability, $(\sigma_c,1]$, with $0<\sigma_c<1$. In this case, the behaviour of the regularized cross-diffusion problem \fer{eq:u1g}-\fer{eq:idg}, for $\delta$ large enough, is similar to that of Turing's classical problem with linear diagonal diffusion \fer{eq:u1:ld}-\fer{eq:u2:ld}.
%\begin{corollary}\label{cor:1}
%	In the conditions of Theorem~\ref{th:estab}, if 
%\begin{align}
%\label{ass:sigmac}
%\tr( K^{-1}D^1(\bu^*))<0 	
%\end{align}
%then there exists $\delta_0>0$ such that for $\delta>\delta_0$ the critical bifurcation parameter satisfies $\sigma_c<1$.
%\end{corollary}
%Condition \fer{ass:sigmac} is also satisfied by competition matrices of the type \fer{def:Bex}. Indeed, for $D^\delta(\bu^*)$ given by \fer{def:ex1} we have $\tr( K^{-1}D^1(\bu^*))=-2/(1-b)$ and for $D^\delta(\bu^*)$ given by \fer{def:ex2} we have
%\begin{align*}
%	 \tr( K^{-1}D^1(\bu^*))=\frac{(b-8)u_1^*u_2^*-2(u_1^*-u_2^*)^2}{4(1-b)u_1^*u_2^*}.
%\end{align*}
%Condition \fer{ass:sigmac} is equivalent to 
%\begin{align*}
% \beta_{22}u_2^* (d_{11}^{11}u_1^*+d_{11}^{21}u_2^*)  &+\beta_{11}u_1^* (d_{22}^{11}u_1^*+d_{22}^{21}u_2^*) \\
% & > \beta_{21}u_2^* (d_{12}^{11}u_1^*+d_{12}^{21}u_2^*) +\beta_{12}u_1^*) (d_{21}^{11}u_1^*+d_{21}^{21}u_2^*).
%\end{align*}

Our second result allows to estimate not only the instability wave numbers provided by the linear analysis but also the amplitude corresponding to these modes. The approximation of the steady state solution is obtained using a weakly nonlinear analysis (WNA) based on the method of multiple scales.
\begin{theorem}
\label{th:wna}
Assume the hypothesis of Theorem~\ref{th:estab} and let $\eps^2=(\delta_c-\delta)/\delta_c$ be a small number. Then, there exist sets of data problem such that the stationary WNA approximation to the solution $\bu$ of problem \fer{eq:u1g}-\fer{eq:idg} is given by  
\begin{align}
\label{wnl}
 \bv(x) =\bu^*+\eps\brho \sqrt{A_\infty} \cos(k_cx)   + \eps^2A_\infty\big(\bv_{20}+
  \bv_{22}\cos(2k_cx)\big) + O(\eps^3),
\end{align}
where $k_c\in\mathbb{Z}$ is the critical wave number corresponding to $\delta_c$, $A_\infty$ is a positive constant and $\brho,\bv_{20}$ and $\bv_{22}$ are constant vectors.
\end{theorem}

Our third result focuses on the limit behaviour of the critical parameters and the amplitude when $\delta\to0$, this is, when the solutions of the approximated problems converge to the solution of the original singular problem. For the sake of simplicity, we limit our study to the following example:
\begin{align}
&  \p_t u_1 - \p_x \big( u_1 (\p_x u_1 +\p_x u_2)\big) = u_1(1-u_1),  \label{eq:u1x2}\\
&  \p_t u_2 -  \p_x \big( u_2 (\p_x u_1 +\p_x u_2)\big) = u_2(4-(2u_1+u_2)), \label{eq:u2x2}
\end{align}
whose solutions we approximate by the two-parameter family of solutions of 
\begin{align}
&  \p_t u_1 - \p_x \big( u_1 ((1+\delta)\p_x u_1 +\p_x u_2)\big) = u_1(1-(u_1+\frac{b}{2}u_2)),  \label{eq:u1x}\\
&  \p_t u_2 -  \p_x \big( u_2 (\p_x u_1 +(1+\delta)\p_x u_2)\big) = u_2(4-(2u_1+u_2)). \label{eq:u2x}
\end{align}
On one hand, Theorem~\ref{th:estab} ensures the existence of $\delta_c>0$ such that, for any $b\geq0$, the equilibrium ${\bu^*=\frac{1}{1-b}(1-2b,2)}$ of \fer{eq:u1x}-\fer{eq:u2x} becomes unstable for $\delta<\delta_c$,  with an associated critical wave number such that $k_c\to\infty$ as $\delta\to0$.

On the other hand, for $\delta<b/4$ and $b\to0$, the sequence of solutions of \fer{eq:u1x}-\fer{eq:u2x}  converges  to a solution of \fer{eq:u1x2}-\fer{eq:u2x2} in the space $BV(0,T,L^\infty(\O))\cup L^\infty(0,T;BV(\O))$.  Therefore, for the approximation \fer{wnl} provided by the weakly nonlinear analysis  to remain  valid for all $\delta>0$, the corresponding amplitude $A_\infty$ must vanish in the limit $\delta\to0$, making in this way compatible the increase of oscillations with its $BV$ regularity. 
 
%
%Finally, we consider the following example to illustrate the behaviour, when $\delta\to0$, of solutions $\bu_\delta$ of problem  \fer{eq:u1g}-\fer{eq:idg} for which $\bu^*$ is unstable:
%\begin{align}
%&  \p_t u_1 - \p_x \big( u_1 ((1+\delta)\p_x u_1 +\p_x u_2)\big) = u_1(1-(u_1+\frac{b}{2}u_2)),  \label{eq:u1x}\\
%&  \p_t u_2 -  \p_x \big( u_2 (\p_x u_1 +(1+\delta)\p_x u_2)\big) = u_2(4-(2u_1+u_2). \label{eq:u2x}
%\end{align}
%On one hand, for $\delta<b/4$ and $b\to0$, the solutions of this sequence of problems, $\bu_\delta$, are expected to converge  to a solution, $\bu$, of 
%\begin{align}
%&  \p_t u_1 - \p_x \big( u_1 (\p_x u_1 +\p_x u_2)\big) = u_1(1-u_1),  \label{eq:u1x2}\\
%&  \p_t u_2 -  \p_x \big( u_2 (\p_x u_1 +\p_x u_2)\big) = u_2(4-(2u_1+u_2). \label{eq:u2x2}
%\end{align}
%in the space $BV(0,T,L^\infty(\O))\cup L^\infty(0,T;BV(\O))$, see \cite{Bertsch2011, Galiano2015}. On the other hand, Theorem~\ref{th:estab} ensures the existence of $\delta_c>0$ such that the equilibrium ${\bu^*=\frac{1}{1-b}(1-2b,2)}$ of \fer{eq:u1x}-\fer{eq:u2x} becomes unstable for $\delta<\delta_c$,  with an associated wave number such that $k_c\to\infty$ as $\delta\to0$.   Therefore, for the approximation provided by the weakly nonlinear analysis to be valid, the corresponding amplitude $A_\infty$ must vanish in the limit $\delta\to0$.

\begin{theorem}
\label{th:example}
	Set $\alpha=(1,4)$, and let $D^\delta(\bu)$ and $B^b$ be  given by \fer{def:ex1} and \fer{def:Bex}, respectively, for $b<1/2$ and  $0<\delta<b/4$. Then, there exists $\delta_c(b)>0$ such that if $\delta<\delta_c(b)$ then 
 $\bu^* =\frac{1}{1-b}(1-2b,2)$ is linearly unstable for problem  \fer{eq:u1g}-\fer{eq:idg}. 
  In addition, 
  \begin{align*}
  \lim_{b\to0} \delta_c(b)=0,\quad 	\lim_{b\to0} k_c(b)=\infty,
  \end{align*}
and the amplitude provided by the weakly nonlinear analysis satisfies
  \begin{align*}
  \lim_{b\to0} A_\infty(b)=0.
  \end{align*}
In particular, the weakly nonlinear approximation $\bv$ given by \fer{wnl} satisfies $\bv\to\bu^*$ uniformly in $\O$ as $b\to0$.
\end{theorem}

\section{Numerical experiments}

In order to analyze the quality of the approximation provided by the WNA, as well as the properties stated in Theorems~\ref{th:estab} to \ref{th:example}, we compare it to a numerical approximation of the evolution problem computed though the finite element method (FEM). 

For the FEM approximation, we used the open source software \texttt{deal.II} \cite{dealII83} to implement a time semi-implicit 
scheme with a spatial linear-wise finite element discretization. For the time discretization, we take 
in the experiments a uniform time partition of time step $\tau=0.01$. For the spatial  discretization, we take a uniform partition of the interval $\O =(0,\pi)$ with spatial step depending on the predicted wave number of the pattern, see Table~\ref{table:1}.

Let, initially, $t=t_0=0$ and  set $(u_1^0, u_2^0)=(u_{10}, u_{20})$. For $n\geq 1$, the discrete problem is: Find $u_1^n, u_2^n \in S^h$ such that  
\begin{align}
 \frac{1}{\tau}\big( u_1^n-u_1^{n-1} , \chi )^h
& + \big( d_{11}^{\delta}(\bu^n) \p_x u_1^n + d_{12}^{\delta}(\bu^n) \p_x u_2^n  ,\p_x\chi \big)^h  =\big(  f_1^b(u_1^n,u_2^n) , \chi )^h, \label{eq:disc1}\\
\frac{1}{\tau}\big( u_2^n-u_2^{n-1} , \chi )^h
& + \big( d_{21}^{\delta}(\bu^n) \p_x u_1^n + d_{22}^{\delta}(\bu^n) \p_x u_2^n  ,\p_x\chi \big)^h  = \big(  f_2^b(u_1^n,u_2^n) , \chi )^h, \label{eq:disc2} 
\end{align}
for every $ \chi\in S^h $, the finite element space of piecewise $\mathbb{Q}_1$-elements. 
Here, $(\cdot,\cdot)^h$ 
stands for a discrete semi-inner product on $\mathcal{C}(\overline{\Omega} )$.

Since \fer{eq:disc1}-\fer{eq:disc2} is a nonlinear algebraic problem, we use a fixed point argument to 
approximate its solution,  $(u_1^n,u_2^n)$, at each time slice $t=t_n$, from the previous
approximation $(u_1^{n-1},u_2^{n-1})$.  Let $u_1^{n,0}=u_1^{n-1}$ and $u_2^{n,0}=u_2^{n-1}$. 
Then, for $k\geq 1$ the linear problem to solve is: Find $(u_1^{n,k},u_2^{n,k})$ such that for for all $\chi \in S^h$ 
\begin{align*}
 \frac{1}{\tau}\big( u_1^{n,k}-u_1^{n-1} , \chi )^h
& + \big( d_{11}^{\delta}(\bu^{n,k-1}) \p_x u_1^{n,k} + d_{12}^{\delta}(\bu^{n,k-1}) \p_x u_2^{n,k}  ,\p_x\chi \big)^h \\
& =\big(  u_1^{n,k} (\alpha_1^b -\beta_{11}^b u_1^{n,k-1} -\beta_{12}^b u_2^{n,k-1}) , \chi )^h, \nonumber\\
\frac{1}{\tau}\big( u_2^{n,k}-u_2^{n-1} , \chi )^h
& + \big( d_{21}^{\delta}(\bu^{n,k-1}) \p_x u_1^{n,k} + d_{22}^{\delta}(\bu^{n,k-1}) \p_x u_2^{n,k}  ,\p_x\chi \big)^h   \\
& =\big(  u_2^{n,k} (\alpha_2^b -\beta_{21}^b u_1^{n,k-1} -\beta_{22}^b u_2^{n,k-1})  \chi )^h. \nonumber
\end{align*}
We use the stopping criteria 
 \begin{equation*}
% \label{tol}
 \max \big(\nor{u_1^{n,k}-u_1^{n,k-1}}_2, \nor{u_2^{n,k}-u_2^{n,k-1}}_2 \big) <\text{tol}_{FP},
 \end{equation*}
for values of $\text{tol}_{FP}$ chosen empirically, and set $(u_1^n, u_2^n)=(u_1^{n,k},u_2^{n,k})$. 
Finally, we integrate in time until a numerical stationary solution, 
$(u_1^S, u_2^S)$, is achieved. This is determined 
by 
\begin{equation*}
%\label{tolS}
\max \big(\nor{u_1^{n,1}-u_1^{n-1}}_2, \nor{u_2^{n,1}-u_2^{n-1}}_2 \big) <\text{tol}_S,
\end{equation*}
where $\text{tol}_S$ is chosen empirically too. In the following experiments we always fix tol$_{FP}=1.e-07$ and tol$_{S}=1.e-12$.

\subsection{Experiment 1}
We investigate the behaviour of the instabilities arising in the solutions of the approximated problems \fer{eq:u1g}-\fer{eq:idg} when $\delta\to0$. Our main aim is to check if the predictions of the weakly nonlinear analysis stated in Theorem~\ref{th:example} are captured by the FEM approximation too. Thus, we use the diffusion matrix $D^\delta(\bu)$ and the competence parameters $B^b$ given by \fer{def:ex1} and \fer{def:Bex}, respectively.

We run three simulations according to the choice of $b$, see Table~\ref{table:1}, and fix  $\delta = 0.95\delta_c(b)$ in all of them, so that $\bu^*$ is unstable and pattern formation follows. 

\begin{table}
\centering
\begin{tabular}{|l|c|c|c|c|}
\hline\hline
   & Simulation 1 & Simulation 2 & Simulation 3 \\ \hline
   $b$& 3.85e-02& 9.91e-03& 4.42e-03\\ \hline
   $\delta(b)$& 4.53e-05& 2.94e-06& 5.83e-07\\ \hline
   $k_c(b)$ & 10& 20& 30\\ \hline
   $A_\infty(b)$ &1.21e-02 & 3.1e-03&1.4e-03 \\ \hline   
   Number of nodes& 128& 256& 512\\ \hline   
 Time steps to stationary& 3.e+04& 1.9e+05& 4.4e+05 \\ \hline          
   Execution time (hours)& 1.67& 19.26& 84.77\\ \hline       
\end{tabular}
 \caption{{\small Data set for the Experiment 1. Wave numbers and times are rounded. Execution time measured for a standard laptop with i7 processor.}}
\label{table:1}
\end{table}

\begin{figure}[h!]
\centering
\includegraphics[width=5cm,height=4cm]{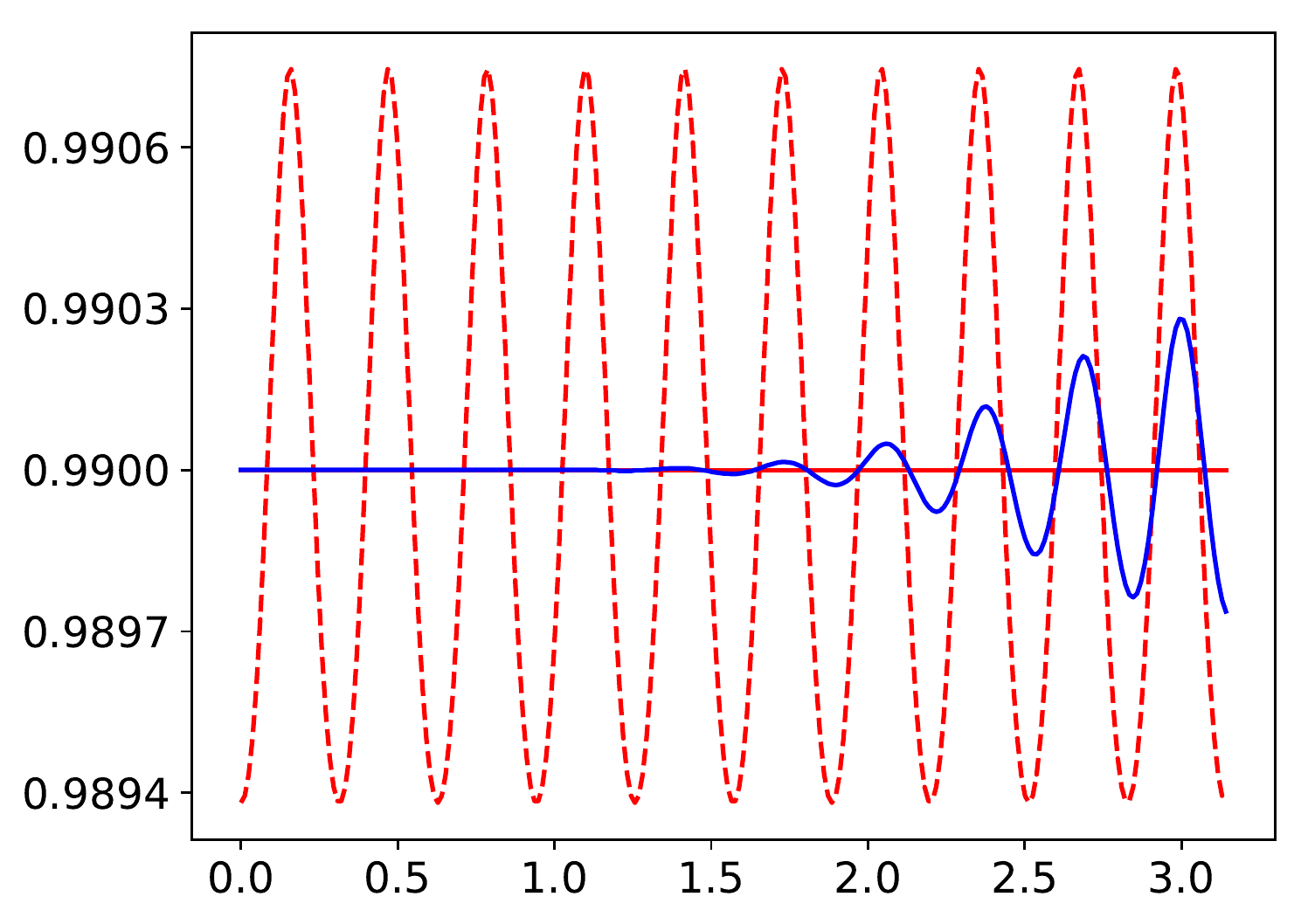}
\includegraphics[width=4.5cm,height=4cm]{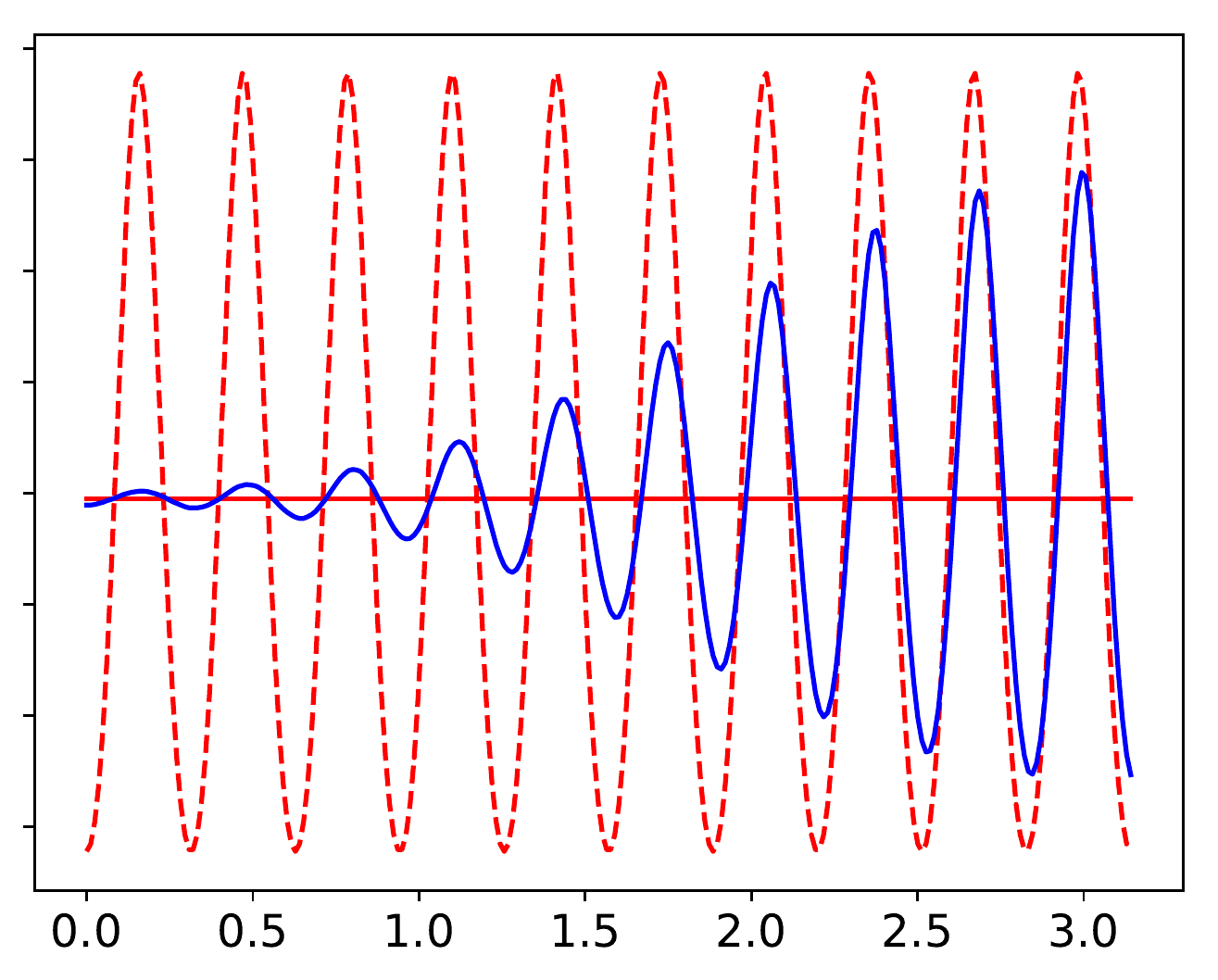}
\includegraphics[width=4.5cm,height=4cm]{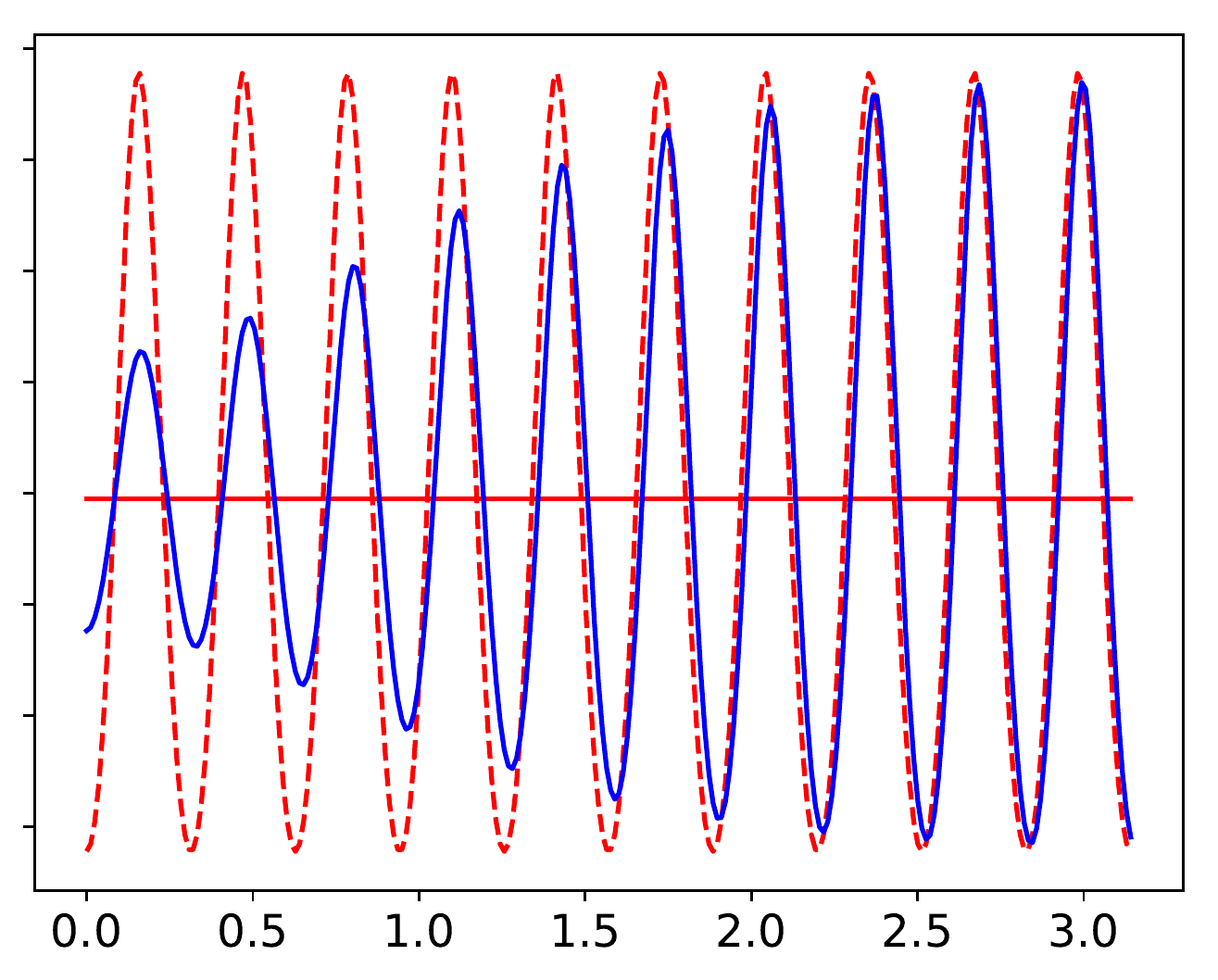}
 \caption{{\small Typical evolution of disturbances 
}} 
\label{fig:1}
\end{figure}

\begin{figure}[h!]
\centering
\includegraphics[width=5cm,height=4cm]{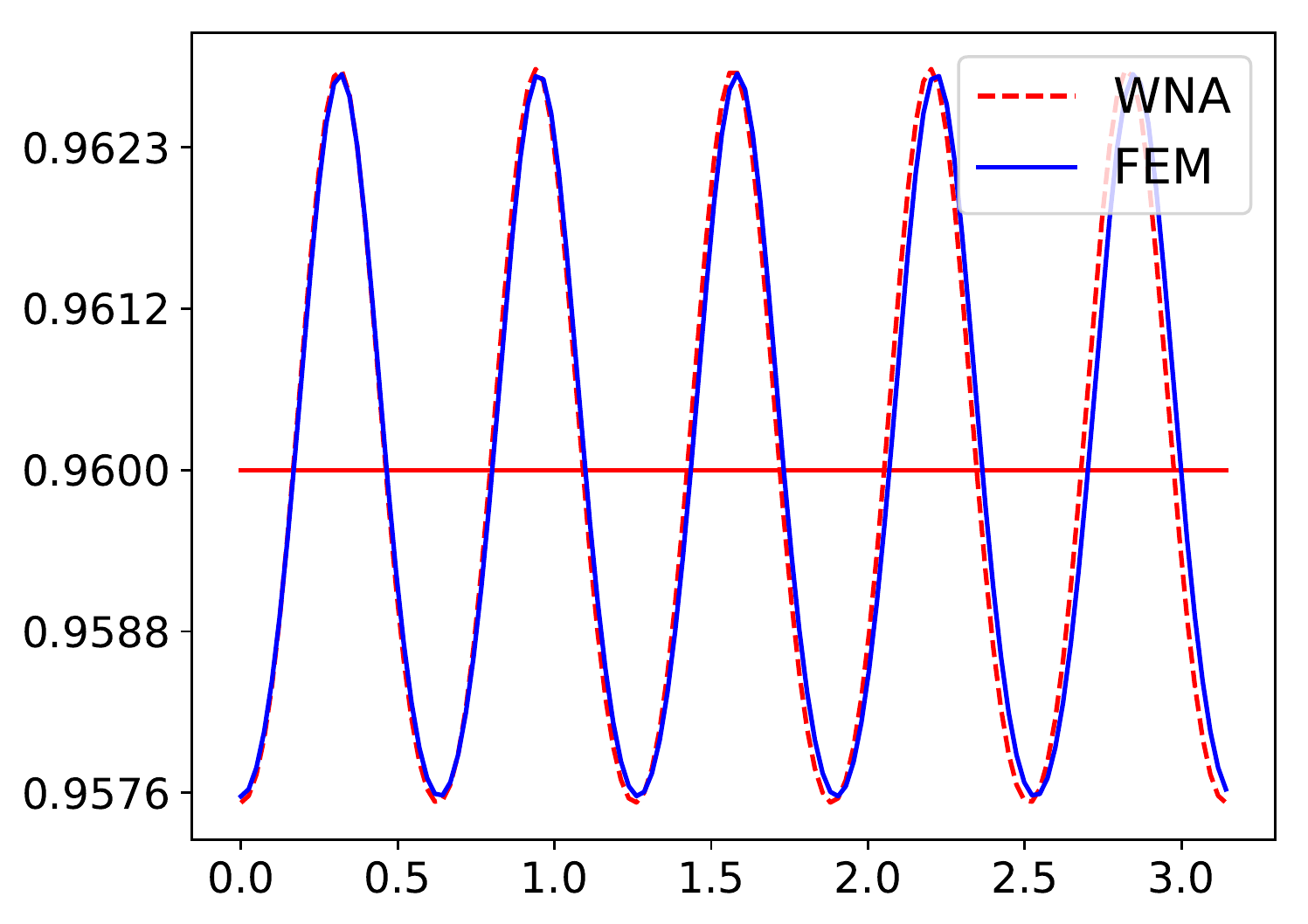}
\includegraphics[width=4.5cm,height=4cm]{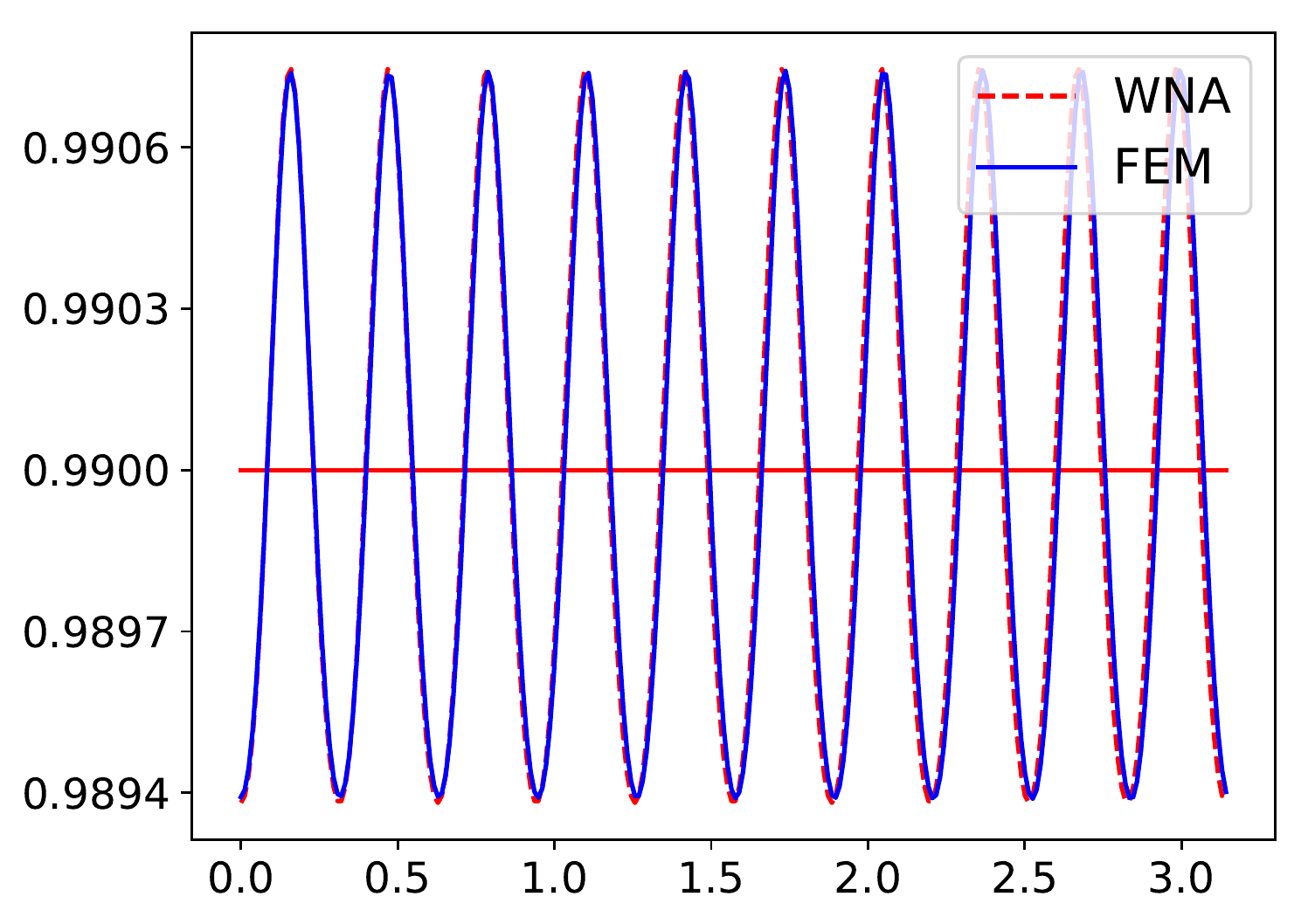}
\includegraphics[width=4.5cm,height=4cm]{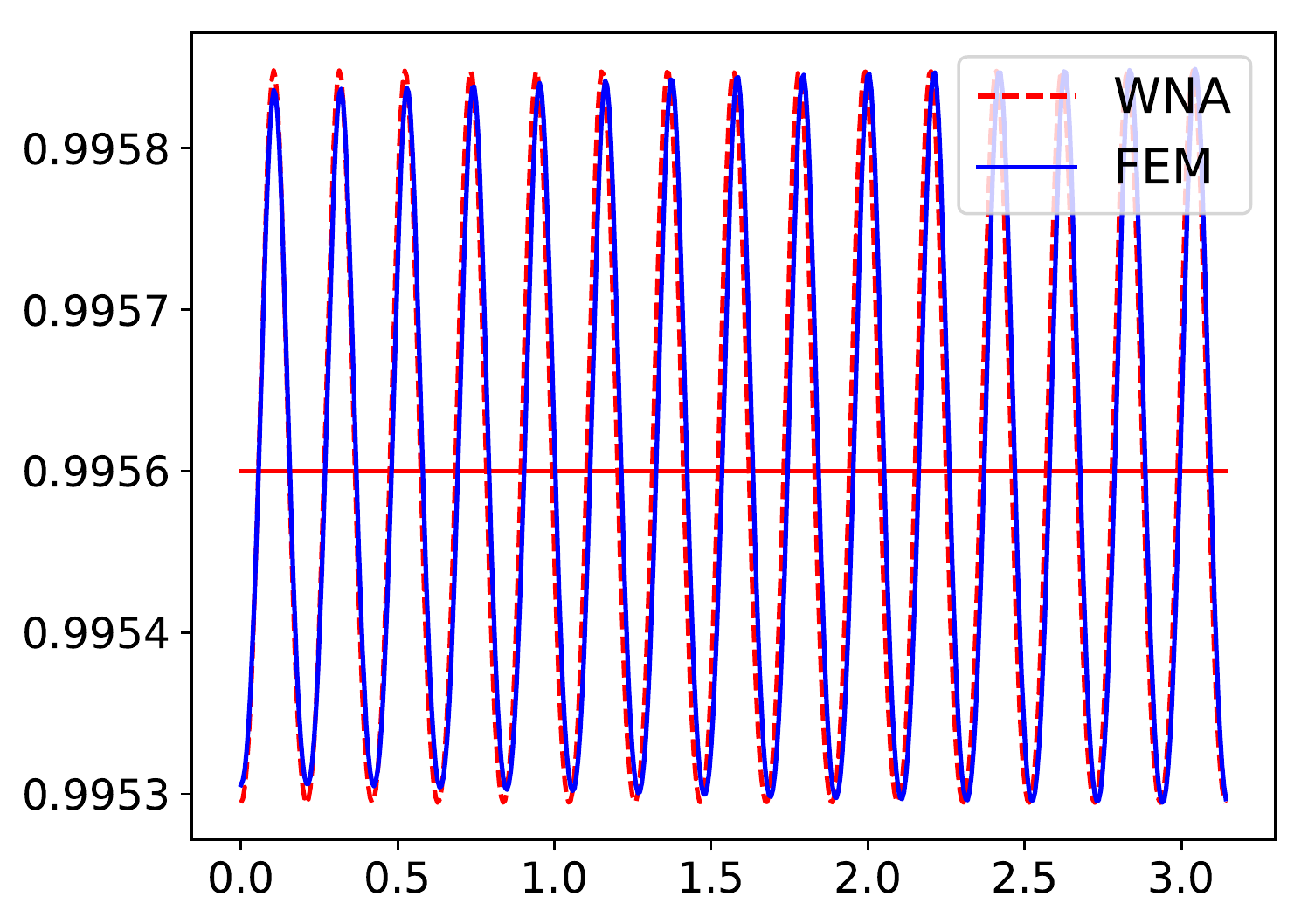}
 \caption{{\small Experiment 1. WNA and FEM approximations corresponding to Simulations 1 to 3 (left to right). Notice the different scales in the ordinates axis showing the decreasing amplitude of the oscillations.}} 
\label{fig:2}
\end{figure}

In Fig.~\ref{fig:1} we show the typical onset and transmission of disturbances found in all the experiments. In this figure and in the following we plot only the first component of the solution, being the behaviour of the second component similar. After a fast decay of the initial data towards the unstable equilibrium, a perturbation with the wave number predicted by the linear analysis grows from one side of the boundary to the rest of the domain until reaching the steady state, see  Fig.~\ref{fig:2}. In the latter figure, we may check the good accordance between the FEM and the WNA approximations which, in numeric figures, have a relative difference of the order $10^{-5}$.

\begin{figure}[h!]
\centering
\includegraphics[width=4.5cm,height=4cm]{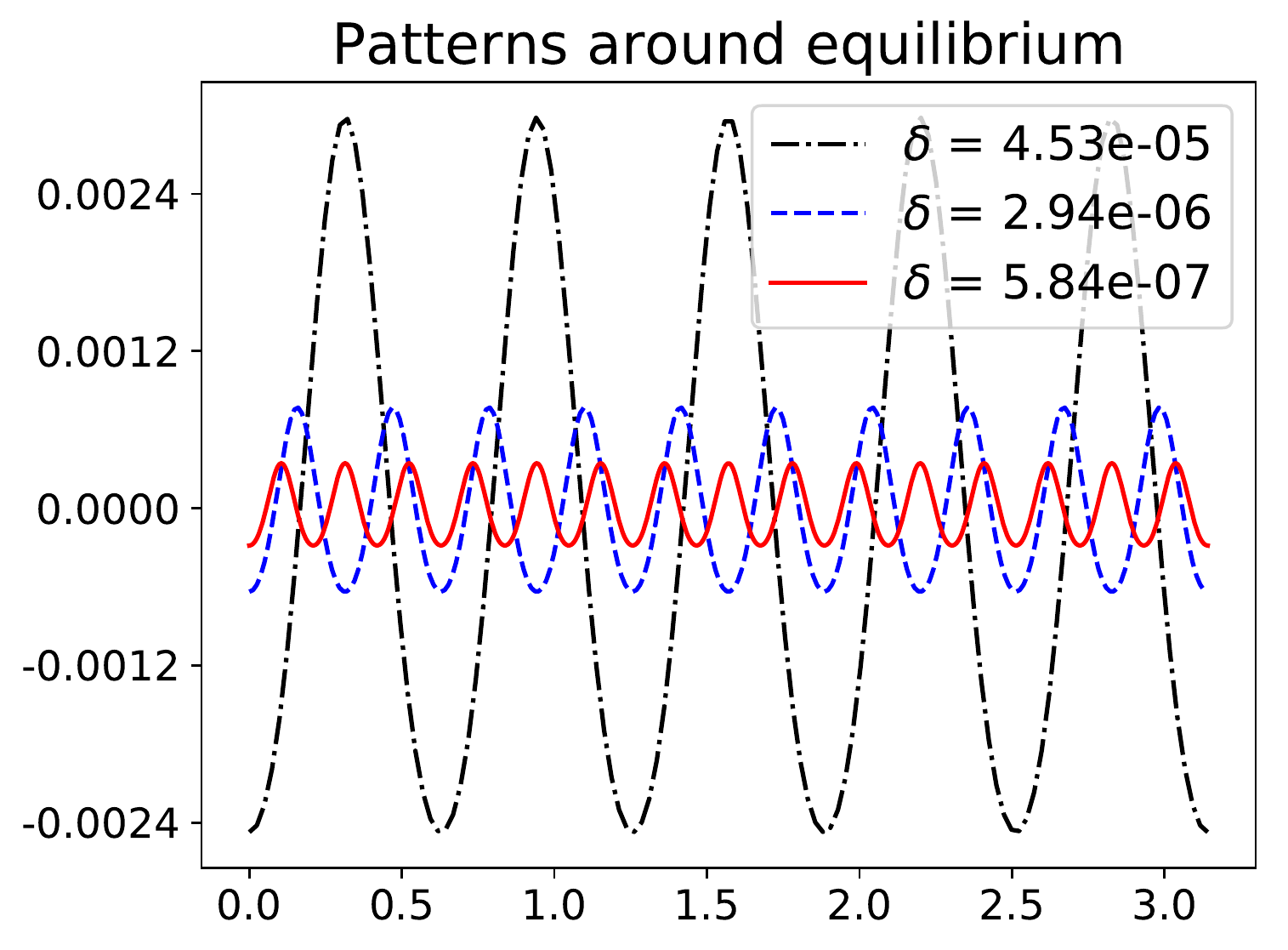}
\includegraphics[width=4.5cm,height=4cm]{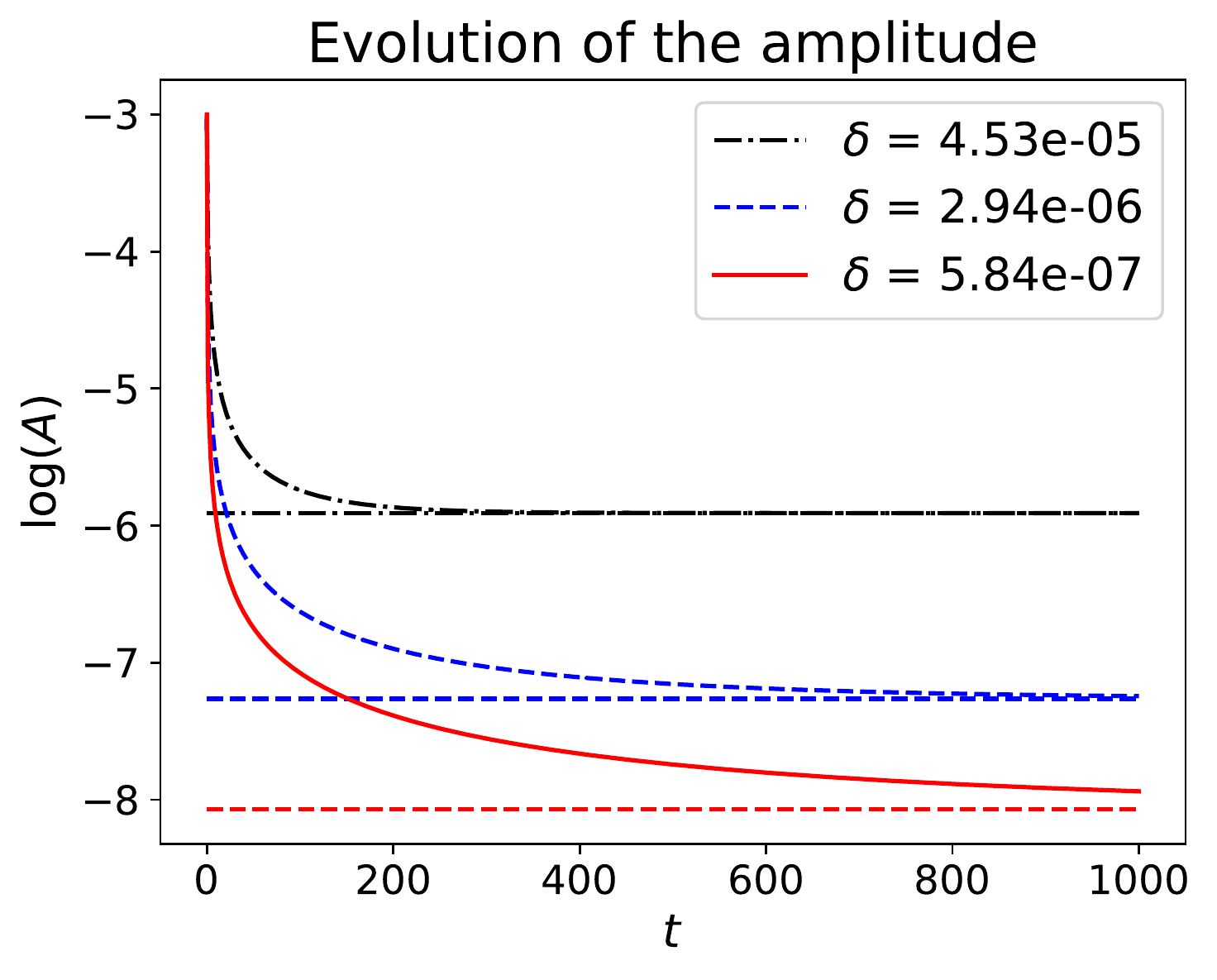}
\includegraphics[width=4.5cm,height=4cm]{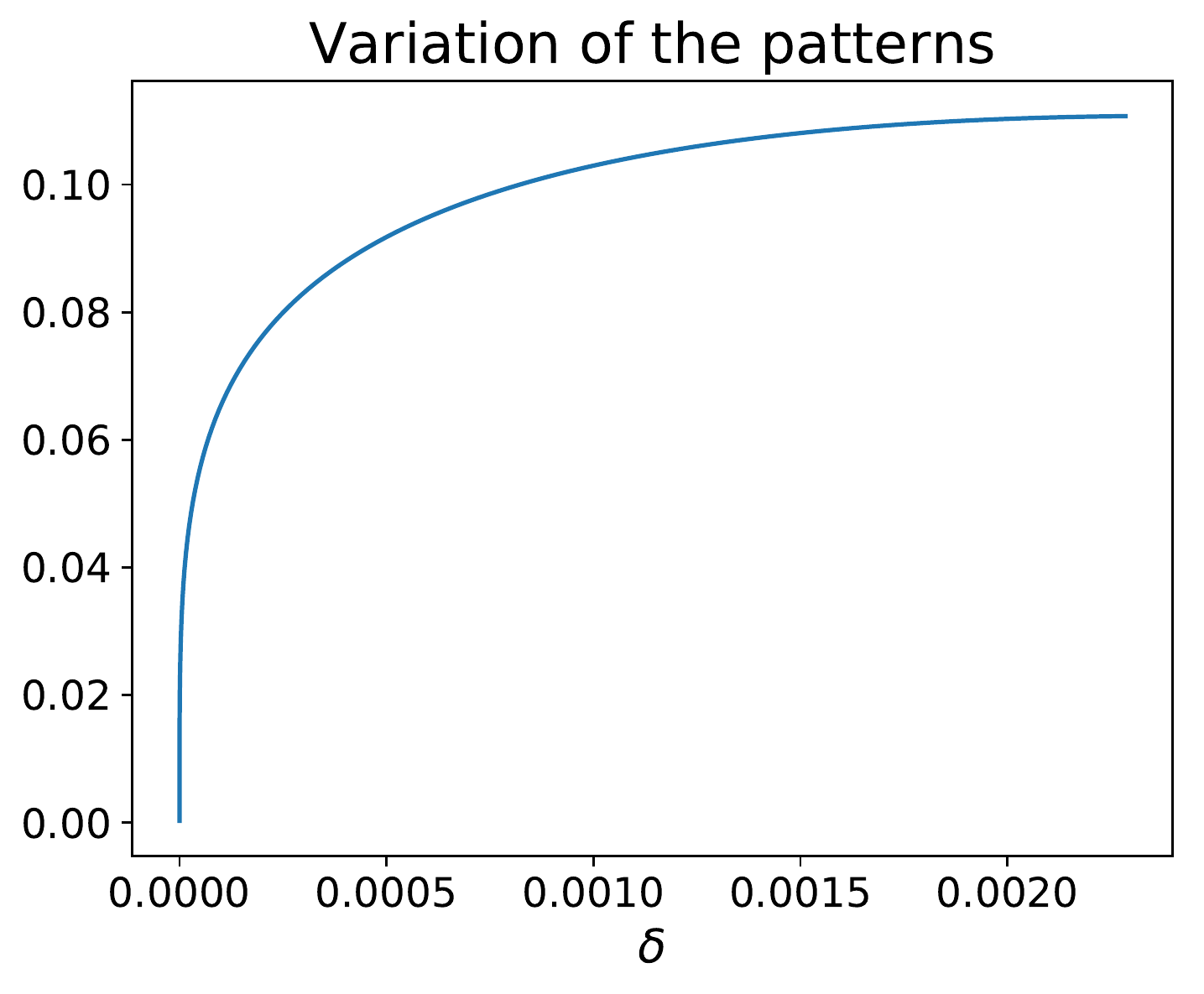}
 \caption{{\small Experiment 1. Behaviour of the patterns as $\delta\to0$.}} 
\label{fig:3}
\end{figure}

In Fig.~\ref{fig:3}  we show three interesting behaviours of solutions when $\delta\to0$. In the left panel, the shrinking amplitude of the stationary patterns while the wave number increases. The equilibrium has been subtracted from the solution to center the pattern in $y=0$. The center panel shows the time evolution of the amplitude (log scale) as given by the exact solution of the Stuart-Landau equation \fer{solex:sl}.  We readily see that the stabilization time is a decreasing function of $\delta$. This fact together with the increment of the wave number when $\delta\to0$ results in very high execution times, see  Table~\ref{table:1}. Finally, the third panel shows how the variation of the numerical stationary solution
\begin{align*}
\int_\O	\abs{\p_x u_1(T,x)}dx
\end{align*}
is an increasing function of $\delta$ and tends to zero as $\delta\to0$, in agreement with the regularity of solutions stated by the theoretical results.

\subsection{Experiment 2}
We repeated Experiment 1 replacing the diffusion matrix $D^\delta(\bu)$ by that defined in  \fer{def:ex2}
In Table~\ref{table:2} we show the relative differences in $L^p$, given by
\begin{align}
\label{reldif}
\RD_p(\vfi_1,\vfi_2) = \frac{\nor{\vfi_1-\vfi_2}_{L^p}}{\nor{\vfi_1}_{L^p}}	,
\end{align}
of the critical bifurcation parameter, $\delta_c$, the \emph{stationary} solution of the FEM  approximation, $\bu(T,\cdot)$, the WNA approximation, $\bv$, and the pattern amplitude, $A_\infty$, corresponding to both approximations of the original diffusion matrix.  We see that although the critical bifurcation parameter is clearly affected by the approximation scheme, the FEM and WNA approximations provided by both schemes are in a very good agreement, as well as the amplitudes of the instability patterns, suggesting that in the limit $\delta\to0$ both sequences of approximations converge to the same limit.

\begin{table}
\centering
\begin{tabular}{|l|c|c|c|c|}
\hline\hline
   & Simulation 1 & Simulation 2 & Simulation 3 \\ \hline
   $\RD_\infty(\delta_c^{(E1)},\delta_c^{(E2)})$& 0.136& 0.117& 0.113 \\ \hline
   $\RD_2(\bu^{(E1)}(T,\cdot),\bu^{(E2)}(T,\cdot))$& 3.74e-06& 8.68e-07& 5.41e-07\\ \hline
   $\RD_2(\bv^{(E1)},\bv^{(E2)})$ & 3.46e-06& 2.13e-07& 4.19e-08\\ \hline
   $\RD_\infty(A_\infty^{(E1)},A_\infty^{(E2)})$ &2.90e-03 & 6.82e-04& 2.99e-4 \\ \hline   
\end{tabular}
 \caption{{\small Comparison between the results obtained with the approximation diffusion matrices corresponding to Example 1 (E1) and Example 2 (E2), given by \fer{def:ex1} and \fer{def:ex2} respectively. RD$_p$ denotes the relative difference in $L^p$, see \fer{reldif}. }}
\label{table:2}
\end{table}

%\vspace{2cm}

%\subsection{Experiment 3}

%\newpage

%%%%%%%%%%%%%%%%%%%%%%%%%%%%%%%%%%%%%%%%%%%%%%%%%%%%%%%%%%%%%%

%\newpage

\section{Proofs}

We use the decomposition of the nonlinear problem \fer{eq:u1g}-\fer{eq:idg}  in terms of its linear and nonlinear parts. Let $\bv=\bu-\bu^*$, where $\bu$ is a solution of \fer{eq:u1g}-\fer{eq:idg}.  Then, $\bv$ satisifies 
\begin{align}
\label{eq:vs}
\p_t \bv = \cL^\delta \bv + \cN^\delta(\bv)	,
\end{align}
where we split the reaction-diffusion terms into their linear parts 
\begin{align*}
%\label{def:ldelta}
 \cL^\delta \bv =   D^\delta(\bu^*)\p_{xx}   \bv +  K\bv , 
\end{align*}
with $K$ given by \fer{def:K}, and their nonlinear parts
\begin{align}
\label{part:nl}
\cN^\delta = \p_x\big(D^\delta(\bv)\p_x\bv\big)	+ \tilde\bf^b(\bv),
\end{align}
being $\tilde f_i^b(\bv)= -\beta_{ii}^bv_i^2-\beta_{ij}^bv_iv_j$, for $i,j=1,2$ and $i\neq j$. 

\bigskip

\no\textbf{Proof of Theorem~\ref{th:estab}. }
We study the linearization of \fer{eq:vs}, this is, the equation 
\begin{align}
\label{eq:linear}
 \p_t \bw = \cL^\delta \bw , 
\end{align}
satisfying Neumann homogeneous boundary conditions and with initial data ${\bw_0=\bu_0-\bu^*}$.
This linear  problem is  well-posed due to the second assumption of $H_D$. The type of boundary conditions lead to seek for solutions of the form 
$ \bw = e^{\lambda t}\cos(kx)\overline{\bw}$, with $k=1,2,\ldots$,
where $\overline{\bw}$ is a constant vector. Replacing $\bw$ in \fer{eq:linear} we obtain the matrix eigenvalue problem 
\begin{align*}
%\label{def:Ak}
 A_k \overline{\bw} = \lambda \overline{\bw}, \qtext{with }A_k = K -k^2  D^\delta(\bu^*).
\end{align*}
Since, by hypothesys, $\tr(A_k)=\tr(K)- k^2\tr(  D^\delta(\bu^*))<0$ for all $k=0,1,\ldots$,  an eigenvalue with positive real part (instability) may exist only if $\det(A_k)$ is negative for some wave number $k$. 
We introduce the notation  $h(k^2) = \det(A_k)$:
\begin{align*}
h(k^2) =  \det(D^\delta(\bu^*)) k^4 + q_\delta (\bu^*) k^2+\det(K), 
%\label{def:detnl}
\end{align*}
where 
$
q_\delta  (\bu^*)  = d_{11}^\delta(\bu^*)\beta_{22}^bu_2^*+ d_{22}^\delta(\bu^*)\beta_{11}^bu_1^*-
 (d_{12}^\delta(\bu^*)\beta_{21}^bu_2^*+ d_{21}^\delta(\bu^*)\beta_{12}^bu_1^*). 
 $
%
%\begin{align}
%  q_\delta  (\bu^*)  = d_{11}^\delta(\bu^*)\beta_{22}u_2^*+ d_{22}^\delta(\bu^*)\beta_{11}u_1^*-
% (d_{12}^\delta(\bu^*)\beta_{21}u_2^*+ d_{21}^\delta(\bu^*)\beta_{12}u_1^*). \nonumber
%\end{align}
%Observe that since $ q_\delta  $ is affine in $\delta$ then it is also monotone with respect to this parameter.
The minimum of the convex parabola $h$ is attained at 
\begin{align*}
%\label{def:km}
 k^2_m(\delta) = -\frac{ q_\delta  (\bu^*)}{2 \det(D^\delta(\bu^*))},
\end{align*}
requiring $ q_\delta  (\bu^*)<0$, which is true in view of \fer{ass:inst2}. A necessary  condition for linear instability is $h(k^2_m(\delta)) <0$, where
\begin{align*}
%\label{eq:min2}
 h(k_m^2(\delta)) = \det(K) - \frac{ q_\delta (\bu^*)^2}{4  \det(D^\delta(\bu^*))}  .
\end{align*}
In this expression, $\det(K)$ is a positive constant and $ q_\delta (\bu^*)^2>0$ for  all $\delta\geq0$. Thus, since $ q_\delta (\bu^*)$ and 
$D^\delta(\bu^*)$ are monotone with respect to $\delta$ and $\det(D^\delta(\bu^*))\to0$ as $\delta\to0$, we deduce the existence of an unique $\bar\delta_c>0$ such that $h(k_m^2(\bar\delta_c)) =0$.% and $h(k_m^2(\delta)) <0$ if $\delta<\bar\delta_c$.  
Therefore,  for $\delta<\bar\delta_c$ we have    $h(k^2(\delta))<0$ if $k^2(\delta) \in (k_-^2(\delta), k_+^2(\delta))$, where
\begin{align*}
%\label{rootsh}
k_{\pm}^2(\delta) = \frac{- q_\delta  (\bu^*)\pm \sqrt{( q_\delta (\bu^*))^2-4 \det(D^\delta(\bu^*))\det(K)}}{2 \det(D^\delta(\bu^*))}.
\end{align*}
Due to the boundary conditions, the onset of instabilities only occurs when one of the extremes values of the interval   $(k_-(\delta), k_+(\delta))$ is an integer number. Since  $k_+(\delta)\to\infty$ as $\delta\to0$, this will certainly holds for $\delta$ small enough. We define 
the \emph{critical bifurcation parameter}, $\delta_c$, as such number, and  the  \emph{critical wave number}, $k_c\in\mathbb{Z}$, as the corresponding root of $h(k^2)$.
Finally, the last assertion of the theorem is a consequence of the infinte limit of $k_+(\delta)$ as $\delta\to0$. $\Box$ 
 %%%%%%%%%%%%%%%%%%%%%%%%%%%%%%%%%%%%%%%%%%%%%%%%%%%%%%%%%%%%%%%%%%%%%%%%%%%%%%%
\bigskip

\no\textbf{Proof of Theorem~\ref{th:wna}. }
We retake the whole nonlinear equation \fer{eq:vs} for $\bv=\bu-\bu^*$. The idea of the weakly nonlinear analysis is to look for an approximation of $\bv$ for a value of $\delta$ near the critical bifurcation parameter $\delta_c$. This approximation is defined as an expansion in terms of a small parameter, that we choose as $\eps^2=(\delta_c-\delta)/\delta_c$, for $\delta<\delta_c$. We consider the expansions
\begin{align*}
& \delta = \delta_c -\eps\delta_1-\eps^2 \delta_2-\eps^3\delta_3+O(\eps^4),\\
& t = \eps t_1+ \eps^2 t_2+ \eps^3t_3+O(\eps^4),\\
& \bv = \eps\bv_1+\eps^2\bv_2+\eps^3\bv_3 + O(\eps^4),	
\end{align*}
and then introduce these expressions in equation \fer{eq:vs} and collect the resulting equations in terms of powers of $\eps$. Since this procedure is standard, we give the results and omit intermediate calculations for the sake of brevity.  We get 
\begin{align}
&\text{Order }\eps:& &\cL^{\delta_c}\bv_1  = 0  .              & \label{prob:e1}\\
&\text{Order }\eps^2:& & \cL^{\delta_c}\bv_2 = \p_{t_1}\bv_1              
              + \delta_1D^1(\bu^*) \p_{xx} \bv_1  & \nonumber \\
  & & &  -\frac{1}{2}\big( \cQ_K(\bv_1,\bv_1) + \p_{xx} \cQ_{D^{\delta_c}}(\bv_1,\bv_1)  \big) - \cS_{D^{\delta_c}} (\bv_1) =:\bF   .        & \label{prob:e2}\\
&\text{Order }\eps^3:& & \cL^{\delta_c}\bv_3 = (\p_{t_1}\bv_2 + \p_{t_2}\bv_1)  
               +D^1(\bu^*)   \p_{xx} (\delta_1\bv_2+\delta_2\bv_1)  -  \cQ_K(\bv_1,\bv_2) \nonumber \\
  & & &  - \p_{xx} \cQ_{D^{\delta_c}}(\bv_1,\bv_2)    
        + \frac{1}{2}\delta_1 \p_{xx}\cR_1(\bv_1) -\cR_2(\bv_1,\bv_2) + \delta_1 \cR_3(\bv_1)  
     =:\bG.     &  \label{prob:e3}
\end{align}
Here,  $D^1(\bu^*)$ is given by \fer{decomp:Ddelta} and 
\begin{align*}
  \cQ_K(\bx,\by) = - \begin{pmatrix}
                            2\beta_{11}^bx_1y_1 +\beta_{12}^b (x_1y_2+x_2y_1) \\
                            2\beta_{22}^bx_2y_2 +\beta_{21}^b (x_1y_2+x_2y_1) 
                           \end{pmatrix},\quad
                           \cQ_{D^\delta}(\bx,\by) = \begin{pmatrix}
                    d_{11}^{\delta  1}x_1y_1 + d_{12}^{\delta  2}x_2y_2\\
                    d_{21}^{\delta  1}x_1y_1 + d_{22}^{\delta  2}x_2y_2
                  \end{pmatrix},
\end{align*}
\begin{align*}                  
& \cS_{D^\delta} (\bv) = \p_x\begin{pmatrix}
                 d_{11}^{\delta  2}v_2 \p_x v_1 + d_{12}^{\delta  1}v_1 \p_x v_2\\
                  d_{21}^{\delta  2}v_2 \p_x v_1 + d_{22}^{\delta  1}v_1 \p_x v_2\\
                 \end{pmatrix}.
\end{align*}
\begin{align*}
\cR_1(\bv_1) = 	\begin{pmatrix}
                    d_{11}^{11}(v_{1 1})^2 + d_{12}^{21}(v_{1 2})^2\\
                    d_{21}^{11}(v_{1 1})^2 + d_{22}^{21}(v_{1 2})^2\end{pmatrix},\quad
\cR_3(\bv_1) = \p_x\begin{pmatrix}
  d_{11}^{21}   v_{12} \p_x v_{11} 
 + d_{12}^{11}   v_{11} \p_x v_{12} 
\\
d_{21}^{21}   v_{12} \p_x v_{11}  
+
 d_{22}^{11}   v_{11} \p_x v_{12} 
 \end{pmatrix},                    
\end{align*}
\begin{align*}
&	\cR_2(\bv_1,\bv_2) = \p_x\begin{pmatrix}
  d_{11}^{\delta_c2}  (v_{12}\p_xv_{21}+v_{22}\p_xv_{11}) 
  +   d_{12}^{\delta_c1}  (v_{11}\p_xv_{22}+v_{21}\p_xv_{12})
\\
  d_{21}^{\delta_c2}  (v_{12}\p_xv_{21}+v_{22}\p_xv_{11}) 
  +   d_{22}^{\delta_c1}  (v_{11}\p_xv_{22}+v_{21}\p_xv_{12})
 \end{pmatrix}, 
\end{align*}
where  we introduced the notation $ d_{ij}^{\delta m} = d_{ij}^{m 0}+\delta d_{ij}^{m 1}$,
%\begin{align}
%\label{def:dijdeltam}
% d_{ij}^{\delta m} = d_{ij}^{m 0}+\delta d_{ij}^{m 1}
%\end{align}
for $m=1,2$ so that  $d_{ij}^\delta(\bv)= d_{ij}^{\delta  1}v_1+d_{ij}^{\delta  2}v_2$. Observe that $d_{ij}^{\delta  1},d_{ij}^{\delta  2}$ are the elements of the matrices $D^{\delta 1},D^{\delta 2}$  introduced in the first assumption of $H_D$. 
Observe also that  \fer{part:nl} may be written as 
\begin{align*}
%\label{decomp:nl}
\cN^\delta \bv =\frac{1}{2}\big(\cQ_K(\bv,\bv) + \p_{xx}\cQ_{D^\delta}(\bv,\bv)\big) + \cS_{D^\delta}(\bv).
\end{align*}
We now compute the solutions corresponding to each order in the expansion. 
%As commented in Remark~\ref{rem:kc}, $k_c$ is not necessarily an integer number and therefore the elements of the Fourier basis corresponding to $k_c$ do not satisfy, in general,  the homogeneous Neumann boundary conditions. We thus replace in our approximated analysis the original spatial domain $(0,\pi)$ by $(0,2\pi/k_c)$, where we shall solve the linear problems \fer{prob:e1}-\fer{prob:e3}.

\no\textbf{Order $\eps$: }The solution of \fer{prob:e1} is given by
\begin{align*}
 \bv_1(t_1,t_2,x) = A(t_1,t_2)\brho \cos(k_cx), \qtext{with }\brho\in \ker( K - k_c^2 D^{\delta_c}(\bu^*)),
\end{align*}
where $A$ is the amplitude of the pattern, unknown at the moment. 
%However, since up to the order $\eps^3$ there only appear the time derivatives $\p_{t_1}$ and $\p_{t_2}$, we do know that $A\equiv A(t_1,t_2)$. 
Observe that  $\ker(A_{k_c}^{\delta_c})$ is a one-dimensional subspace, implying that the vector $\brho$ is defined up to a multiplicative constant. We shall fix this constant later. 

\medskip

\no\textbf{Order $\eps^2$: }We start expressing $\bF$ in terms of $A$ and $\brho$.
%\begin{align*}
%\bF = \p_{t_1}\bv_1              
%              + \delta_1D^1(\bu^*) \p_{xx} \bv_1    -\frac{1}{2}\big(\gamma \cQ_K(\bv_1,\bv_1) + \p_{xx} \cQ_{D^{\delta_c}}(\bv_1,\bv_1)  \big) - \cS_{D^{\delta_c}} (\bv_1) .       
%\end{align*}
We have 
\begin{align*}
& \p_{t_1}\bv_1 = \p_{t_1}A \cos(k_cx) \brho\\
&  \delta_1 D^1(\bu^*)   \p_{xx} \bv_1 = 
 - \delta_1 A  k_c^2 \cos(k_c x) D^1(\bu^*)\brho
\end{align*}
On noting that $\cQ_U(\bv_1,\bv_1) = A^2 \cQ_U(\brho,\brho) \cos^2(k_cx)$, for $U=K,~D^{\delta_c}$, we find 
\begin{align*}
 \frac{1}{2}\big(\cQ_K(\bv_1,\bv_1) +  \p_{xx} \cQ_{D^{\delta_c}}(\bv_1,\bv_1) \big) =&  \frac{1}{4}A^2 \sum_{j=0,2}\cM_j(\brho,\brho) \cos(jk_cx),
\end{align*}
with $ \cM_j =  \cQ_K - j^2k_c^2 \cQ_{D^{\delta_c}}$.
%\begin{remark}
%For $\bu = (u_1,u_2)=(\xi_1 \cos(\alpha_1 x), \xi_2 \cos(\alpha_2 x))$, we have
%\begin{align}
%\label{def:formula}
% \p_x(u_2\p_xu_1) = -\frac{1}{2}\xi_1\xi_2\alpha_1 \big((\alpha_1+\alpha_2) \cos((\alpha_1+\alpha_2)x) + 
% (\alpha_1-\alpha_2) \cos((\alpha_1-\alpha_2)x) \big).
%\end{align}
%\end{remark}
Using standard trigonometric identities, we get 
$ \cS_{D^{\delta_c}} (\bv_1) = - k_c^2 A^2\rho_1\rho_2 \cos(2k_cx)\bd$,
where $\bd = (d_{11}^{\delta_c  2} + d_{12}^{\delta_c  1},  
d_{21}^{\delta_c  2} + d_{22}^{\delta_c  1})$. 
Gathering the above expressions, we obtain 
\begin{align*}
 \bF = &    \Big[   \p_{t_1} A \brho - \delta_1 A  k_c^2 D^1(\bu^*)    \brho \Big] \cos(k_cx) 
  - \frac{1}{4}A^2 \sum_{j=0,2}\cM_j(\brho,\brho) \cos(jk_cx) \\
&   + k_c^2 A^2\rho_1\rho_2 \bd \cos(2k_cx).
\end{align*}
By Fredholm's alternative, \fer{prob:e2} admits a solution if and only if 
$\langle\bF,\bpsi\rangle_{L^2} = 0$, where $\langle\cdot,\cdot\rangle_{L^2} $ denotes the scalar product in $L^2(0,\pi)$, and $\bpsi \in \ker( (\cL^{\delta_c})^*)$ is of the form 
\begin{align}
\label{def:psi}
\bpsi = \boldeta \cos(k_cx), \qtext{with }\boldeta \in \ker( ( K - k_c^2 D^{\delta_c}(\bu^*))^*). 
\end{align}
Observe that $\boldeta$, for similar reasons than $\brho$, is defined up to a multiplicative constant. We fix $\boldeta$ at the end of this proof, and also show that $\langle \brho,\boldeta\rangle \neq 0$.

The compatibility condition implies 
%\begin{align*}
% 0=\langle\bF,\bpsi\rangle_{L^2} = & ~  \langle  (\p_{t_1} A \brho- \delta_1 A  k_c^2 D^1(\bu^*)    \brho  ,\boldeta\rangle    \int_0^{\frac{2\pi}{k_c}} \cos^2 (k_cx) dx   \\
% - & \frac{1}{4} A^2 \sum_{j=0,2} \langle   \cM_j(\brho,\brho)  , \boldeta \rangle ~ \int_0^{\frac{2\pi}{k_c}} \cos(jk_cx) \cos(k_cx) dx\\
%  + &  k_c^2 A^2\rho_1\rho_2  \langle \bd , \boldeta\rangle \int_0^{\frac{2\pi}{k_c}} \cos(2k_cx) \cos(k_cx) dx ~ . 
%\end{align*}
%The second and third integrals at the right hand side vanish, so we obtain
\begin{align*}
    \p_{t_1} A(t_1,t_2) = \delta_1 k_c^2 \frac{\langle D^1(\bu^*)\brho ,\boldeta\rangle}{\langle \brho   ,\boldeta\rangle} A(t_1,t_2).
\end{align*}
Since the solution to this equation is an exponential function, we do not obtain from it any useful indication on the asymptotic behaviour of the pattern amplitude. Therefore, to suppress the secular terms appearing in $\bF$, we impose
\begin{align}
\label{ass:order1}
t_1\equiv 0 \qtextq{and} \delta_1\equiv 0	.
\end{align}
In particular, this implies $A\equiv A(t_2)$.

Assuming these restrictions, the Fredholm's alternative is satisfied, and motivated by the functional form of $\bF$, we seek for a solution of \fer{prob:e2} of the form
\begin{align*}
 \bv_2(t_2,x) = A^2(t_2)\sum_{j=0,2} \bv_{2j}\cos(jk_cx),
\end{align*}
where $\bv_{2j}$ are constant vectors. The linear operator $\cL^{\delta_c}$ may be decomposed as
\begin{align*}
 \cL^{\delta_c} \bv_2=  A^2\sum_{j=0,2} L_j\bv_{2j}\cos(jk_cx), \qtext{with }L_j =  K- j^2k_c^2D^{\delta_c}(\bu^*).
\end{align*}
Then, $\cL^{\delta_c} \bv_2 =\bF$ if the vectors $\bv_{2j}$ are the solutions of the linear systems
\begin{align*}
%\label{def:ws}
 &L_0 \bv_{20} = -\frac{1}{4} \cM_0(\brho,\brho),\quad 
 L_2 \bv_{22} = k_c^2  \rho_1\rho_2  \bd -\frac{1}{4} \cM_2(\brho,\brho).
\end{align*}

\medskip

\no\textbf{Order $\eps^3$: }
We have to solve $\cL^{\delta_c}\bv_3 = \bG$, where, taking into account \fer{ass:order1}, 
\begin{align*}
\bG =&  \p_{t_2}\bv_1  
               +\delta_2D^1(\bu^*)   \p_{xx} \bv_1  -  \cQ_K(\bv_1,\bv_2) 
    - \p_{xx} \cQ_{D^{\delta_c}}(\bv_1,\bv_2)    
         -\cR_2(\bv_1,\bv_2).   	
\end{align*}
Replacing the solutions obtained for the orders $\eps$ and $\eps^2$, i.e.  
$\bv_1 = A(t_2) \brho \cos(k_cx)$ and $\bv_2 = A(t_2)^2 ( \bv_{20} + \bv_{22}\cos(2k_cx))$
in  $\bG$ yields
\begin{align*}
 \bG =&  \Big(\brho \p_{t_2} A  -A k_c^2 \delta_2D^1(\bu^*)\brho -A^3  \big(\cM_1(\brho,\bv_{20}) +  \frac{1}{2} \cM_1(\brho,\bv_{22}) +k_c^2\bR_1\big)\Big) \cos(k_c x)\\
 & - A^3 \Big( \frac{1}{2} \cM_3(\brho,\bv_{22})+ k_c^2\bR_2 \Big)  \cos(3k_cx),
\end{align*}
where
\begin{align*}
\bR_1^{(i)} = &  	d_{i1}^{\delta_c2} \Big[\rho_1 \Big(\frac{1}{2}\bv_{22}^{(2)} -\bv_{20}^{(2)}\Big)-\rho_2 \bv_{22}^{(1)} \Big]
+ d_{i2}^{\delta_c1}  \Big[\rho_2 \Big(\frac{1}{2}\bv_{22}^{(1)} -\bv_{20}^{(1)}\Big)  - \rho_1 \bv_{22}^{(2)}\Big] ,\\
\bR_2^{(i)} = & -3 \Big( d_{i1}^{\delta_c2}\Big[\rho_2 \bv_{22}^{(1)} +
\frac{1}{2}\rho_1 \bv_{22}^{(2)} \Big]+  d_{i2}^{\delta_c1}  \Big[ \rho_1 \bv_{22}^{(2)} + \frac{1}{2}\rho_2 \bv_{22}^{(1)}\Big]\Big) .
\end{align*}

The solvability condition for problem \fer{prob:e3} is $\langle\bG,\bpsi\rangle_{L^2} = 0$, with $\bpsi =\boldeta\cos(k_cx)$ given by \fer{def:psi}. This condition leads to the differential equation
\begin{align*}
  \langle \brho,\boldeta\rangle \p_{t_2} A = \langle\bG_1,\boldeta\rangle A +
  \langle\bG_2 ,\boldeta\rangle A^3 ,
\end{align*}
where
\begin{align}
 \bG_1 = & k_c^2 \delta_2D^1(\bu^*)\brho , \label{def:G1}\\
 \bG_2 =& \cM_1(\brho,\bv_{20}) +  \frac{1}{2} \cM_1(\brho,\bv_{22}) +k_c^2\bR_1,\nonumber %\label{def:G2}
\end{align}
Thus, we deduce the cubic Stuart-Landau equation for the amplitude
\begin{align}
\label{eq:sl}
 \p_{t_2} A = \sigma A - \ell A^3,
\end{align}
with
% \begin{align*}
%  \sigma = \frac{\langle\bG_1^{(1)},\boldeta\rangle}{\langle\brho,\boldeta\rangle}, 
%  \qquad
%  L =   \frac{\langle\bG_3^{(1)},\boldeta\rangle}{\langle\brho,\boldeta\rangle},
% \end{align*}
\begin{align}
\label{def:sL}
 \sigma = \frac{\langle \bG_1, \boldeta\rangle}{\langle \brho ,\boldeta\rangle}, 
 \qquad
 \ell =  - \frac{\langle \bG_2, \boldeta\rangle}{\langle \brho ,\boldeta\rangle}. 
\end{align}

We, finally, fix the vectors $\brho\in \ker( K - k_c^2 D^{\delta_c}(\bu^*))$ , and 
$\boldeta \in \ker( ( K - k_c^2 D^{\delta_c}(\bu^*))^*)$. Since all the elements of both matrices are negative, we may set $ \brho = (1,M)^t$ and  $\boldeta =(1,M^*)^t$ for some $M,M^*<0$, implying $ \langle \brho ,\boldeta\rangle>0$. Thus, the asymptotic behaviour of the solution to \fer{eq:sl} is fully determined by the signs of the numerators in \fer{def:sL}. 

When $\sigma$ and $\ell$ are positive, the amplitude estabilizes to a positive value, this is, $A(t_2)\to A_\infty:=\sqrt{\sigma/\ell}$ as $t_2\to\infty$.
Therefore, in this case, the corresponding solution $ \bv = \eps \bv_1 + \eps^2 \bv_2 + O(\eps^3)$, is given by  
\begin{align*}
 \bv =\eps\brho \sqrt{\frac{\sigma}{\ell}} \cos(k_cx) 
  + \eps^2\frac{\sigma}{\ell}\big(\bv_{20}+
  \bv_{22}\cos(2k_cx)\big) + O(\eps^3). 
\end{align*}
An example of this situation is studied in Theorem~\ref{th:example}. $\Box$
 
%%%%%%%%%%%%%%%%%%%%%%%%%%%%%%%%%%%%%%%%%%%%%%%%%%%%%%%%%%%%%%%%%%%%%%%%%%%%%%%%%%%%%%%%%%%%
\bigskip

\no\textbf{Proof of Theorem~\ref{th:example}. }
%For the specific data, we have 
%\begin{align*}
%\bu^*=\frac{1}{1-b}
%\begin{pmatrix}
%	1-2b\\ 2
%\end{pmatrix},\quad
%K = \begin{pmatrix}
% -u_1^* & -\frac{b}{2}u_1^*\\ -2u_2^*& -u_2^*	
% \end{pmatrix},\quad
%K^{-1} = \frac{1}{u_1^* u_2^* (1-b)}\begin{pmatrix}
% -u_2^* & \frac{b}{2}u_1^*\\ 2u_2^*& -u_1^*	
% \end{pmatrix}.
%\end{align*}
Our aim is to compute the coefficients of the Stuart-Landau equation \fer{eq:sl}. Specifically, we are interested in the ratio  
\begin{align*}
%\label{def:sl}
\frac{\sigma}{\ell} = -\frac{\langle \bG_1, \boldeta\rangle}{\langle \bG_2, \boldeta\rangle}.
\end{align*}

%\bigskip

\no\textbf{Determination of $\langle \bG_1, \boldeta\rangle$.}
\no For the given data, we get $q_\delta(\bu^*) = -u_1^*u_2^* \Big(\frac{b}{2}-2\delta\Big)$,
which is negative if $b>4\delta$. The corresponding roots of $h(k_m^2)$ are positive
and, therefore, we take $\delta_c = \delta_-$, so that for any $\delta<\delta_c$ we have $h(k_m^2)<0$. 
The corresponding critical wave number is the minimum of $h(k^2)$, given by 
\begin{align*}
k_c^2 = 	\frac{b-4\delta_c}{4\delta_c(2+\delta_c)}.
\end{align*}
The vectors $\brho=(1,M)$ and $\boldeta = (1,M^*)$ are elements of $\ker(A_{k_c}^{\delta_c})$ and $\ker(A_{k_c}^{\delta_c})^*$, respectively. Thus, 
\begin{align*}
M=-\frac{1+k_c^2(1+\delta_c)}{\frac{b}{2}+k_c^2}\quad 	M^*=-\frac{u_1^*(1+k_c^2(1+\delta_c))}{u_2^*(2+k_c^2)}.
\end{align*}
From \fer{def:G1}, we obtain $\bG_1 = k_c^2 \delta_2 (u_1^*,u_2^*M)$,
and then $\langle \bG_1, \boldeta\rangle = 	k_c^2 \delta_2 (u_1^*+u_2^*MM^*)$.

\begin{lemma}\label{lemma:orders}
Let $\eps_M=1+M$. We have:
\begin{align}
&	\lim_{b\to0}\delta_c(b) = \lim_{b\to0}\delta_c'(b) =\lim_{b\to0}\delta_c k_c^2 =0,\quad
 \lim_{b\to0}k_c^2 = \infty, \label{lemma:orders:1}\\
& \lim_{b\to0} M = -1,\quad 	\lim_{b\to0} M^* =-\frac{1}{2}, \label{lemma:orders:2}\\
& \lim_{b\to0} k_c^2\eps_M = -1,\quad \lim_{b\to0} k_c^2(1+2M^*)= 9.  \label{lemma:orders:3} 
\end{align}
\end{lemma}
Taking into account that $\delta_2\approx \delta_c$, a first consequence of  Lemma~\ref{lemma:orders} is 
\begin{align}
\label{lim:G1}
\lim_{b\to0} 	\langle \bG_1, \boldeta\rangle = 0.
\end{align}

%%%%%%%%%%%%%%%%%%%%%%%%%%%%%%%%%%%%%%%%%%%%%%%%%%%%%%%%%%%%%%%%%%%%%%%%%%%%%%%%%%%%
\bigskip

\no\textbf{Proof of Lemma~\ref{lemma:orders}. } For proving \fer{lemma:orders:1}, 
we use L'H\^opital's rule to get 
\begin{align*}
\lim_{b\to0}\delta_c(b) = \frac{1}{4}	\lim_{b\to0}\Big(-3-\frac{-6(4-3b)-3b^2}{2\sqrt{(4-3b)^2-b^3}}\Big)=0.
\end{align*}
Let
$
\vfi(b,\delta)=4\det(D^\delta(\bu^*)\det(K)-q_\delta^2(\bu^*)=-4b\delta^2+(8-6b)\delta-b^2/4	
$.
By definition of $\delta_c$, we have $\vfi(b,\delta_c(b))=0$ for all $b\in(0,1)$. Thus
 \begin{align}
 \label{der:dc}
 0=\dfrac{d}{db}\vfi(b,\delta_c(b))=\p_b \vfi(b,\delta_c(b)) + \p_\delta \vfi(b,\delta_c(b))\delta_c'(b).	
 \end{align}
Since $\p_b \vfi(0,0)=0$ and  $\p_\delta \vfi(0,0)=8$, we deduce $\delta_c'(0)=0$. 
We then have
\begin{align*}
\lim_{b\to0} k_c^2 = 	\lim_{b\to0} \frac{b-4\delta_c}{4\delta_c(2+\delta_c)} = \lim_{b\to0}
\frac{1-4\delta_c'}{8\delta_c'(1+\delta_c))}=\infty,
\quad \lim_{b\to0}\delta_c k_c^2 = 	\lim_{b\to0}\frac{b-4\delta_c}{4(2+\delta_c)} = 0. 
\end{align*}
%and
%\begin{align*}
%\lim_{b\to0}\delta_c k_c^2 = 	\lim_{b\to0}\frac{b-4\delta_c}{4(2+\delta_c)} = 0. 
%\end{align*}
The limits \fer{lemma:orders:2} follow easily from the definitions of $M$ and $M^*$. 
Finally, for proving \fer{lemma:orders:3}, we use the defintion of $M$ to get 
%
%\begin{align*}
%M=-\frac{1+k_c^2(1+\delta_c)}{\frac{b}{2}+k_c^2}\quad 	M^*=-\frac{u_1^*(1+k_c^2(1+\delta_c))}{u_2^*(2+k_c^2)}.
%\end{align*}
\begin{align*}
k_c^2\eps_M = 	k_c^2\frac{\frac{b}{2}+k_c^2-1-k_c^2(1+\delta_c)}{\frac{b}{2}+k_c^2} = 
\frac{\frac{b}{2}-1-k_c^2\delta_c}{\frac{b}{2k_c^2}+1},
\end{align*}
from where the first limit follows. The second limit is computed in a similar way. We write
\begin{align*}
k_c^2(1+2M^*) = k_c^2\Big(1+\frac{u_2^*}{u_1^*}M^*	+ M^*\Big(2-\frac{u_2^*}{u_1^*}\Big)\Big).
\end{align*}
On one hand, we have 
\begin{align*}
k_c^2\Big(1+\frac{u_2^*}{u_1^*}M^*\Big) = k_c^2\frac{2+k_c^2-(1+k_c^2(1+\delta_c))}{2+k_c^2}
= k_c^2\frac{1-k_c^2\delta_c}{2+k_c^2}\to 1 \qtext{as }b\to0.
\end{align*}
On the other hand, using the definition of $\bu^*$ and $k_c^2$, we obtain
\begin{align*}
k_c^2\Big(2-\frac{u_2^*}{u_1^*}\Big) =  -\frac{1}{1-2b}4k_c^2b,\qtextq{with}
4k_c^2b = \frac{1}{2+\delta_c}\Big(\frac{b^2}{\delta_c}-4b\Big).
\end{align*}
Using a concatenation of L'H\^opital's rule, we get $\lim_{b\to0}b^2/\delta_c(b) = 	2\lim_{b\to0}1/\delta_c''(0).$
%\begin{align*}
%\lim_{b\to0}\frac{b^2}{\delta_c(b)} = 	2\lim_{b\to0}\frac{1}{\delta_c''(0)}.
%\end{align*}
Differentiating \fer{der:dc} with respect to $b$ yields $\delta_c''(0)=1/16$, implying the result. $\Box$
%%%%%%%%%%%%%%%%%%%%%%%%%%%%%%%%%%%%%%%%%%%%%%%%%%%%%%%%%%%%%%%%%%%%%%%%%%%%%%%%%%%%

\medskip

\no\textbf{Determination of $\langle \bG_2, \boldeta\rangle$.}
The following lemma gives the expression of this scalar product. Since the calculation is  straightforward, we omit the details.
\begin{lemma}
\label{lemma:G2}
We have 
	\begin{align*}
-\langle \bG_2,\boldeta \rangle =S_1 \bv_{20}^{(1)}+S_2 \bv_{20}^{(2)}+T_1  \bv_{22}^{(1)}+T_2   \bv_{22}^{(2)},
\end{align*}
where
\begin{align*}
	\bv_{20} & = \frac{1}{4u_1^*u_2^*(1-b)}\begin{pmatrix}
-u_2^*(2+bM)+bu_1^*M(M+2)\\
2u_2^*(2+bM)-2u_1^*M(M+2)
\end{pmatrix},\\
\bv_{22} & = \frac{1}{\det(L_2)}\big(\bw_2 k_c^4+\bw_1 k_c^2+\bw_0	\big),
\end{align*}
with
\begin{align*}
\bw_2 =  &   4\begin{pmatrix}
 -(1+\delta_c)u_2^*	(\eps_M+\delta_c) + u_1^*M(\eps_M+M\delta_c)\\
u_2^*(\eps_M+\delta_c) - (1+\delta_c)u_1^*M(\eps_M+M\delta_c)
\end{pmatrix},\\
\bw_1 = & \begin{pmatrix}
 -4(1+\delta_c)u_2^*	(\frac{1}{2}+\frac{bM}{4}) -u_2^*(\eps_M+\delta_c)+4u_1^*(M+\frac{M^2}{2})+\frac{b}{2}u_1^*M(\eps_M+M\delta_c)\\
4u_2^*(\frac{1}{2}+\frac{bM}{4})+2u_2^*(\eps_M+\delta_c)-4(1+\delta_c)u_1^*(M+\frac{M^2}{2})-u_1^*M(\eps_M+M\delta_c)
\end{pmatrix},\\
\bw_0= &  \begin{pmatrix}
-u_2^*(\frac{1}{2}+\frac{bM}{4})+\frac{b}{2}u_1^*(M+\frac{M^2}{2})\\
2u_2^*(\frac{1}{2}+\frac{bM}{4})-u_1^*(M+\frac{M^2}{2})
\end{pmatrix},
\end{align*}
%\begin{align*}
%	\bv_{22} 
%& =  \frac{4k_c^4}{\det(L_2)}\begin{pmatrix}
% -(1+\delta_c)u_2^*	(1+M+\delta_c) + u_1^*(M+M^2(1+\delta_c)\\
%u_2^*(1+M+\delta_c) - (1+\delta_c)u_1^*(M+M^2(1+\delta_c)
%\end{pmatrix}\\
%& +  \frac{k_c^2}{\det(L_2)}\begin{pmatrix}
% -4(1+\delta_c)u_2^*	(\frac{1}{2}+\frac{bM}{4}) -u_2^*(1+M+\delta_c)+4u_1^*(M+\frac{M^2}{2})+\frac{b}{2}u_1^*(M+M^2(1+\delta_c)\\
%4u_2^*(\frac{1}{2}+\frac{bM}{4})+2u_2^*(1+M+\delta_c)-4(1+\delta_c)u_1^*(M+\frac{M^2}{2})-u_1^*(M+M^2(1+\delta_c)
%\end{pmatrix}\\
%& +  \frac{1}{\det(L_2)}\begin{pmatrix}
%-u_2^*(\frac{1}{2}+\frac{bM}{4})+\frac{b}{2}u_1^*(M+\frac{M^2}{2})\\
%2u_2^*(\frac{1}{2}+\frac{bM}{4})-u_1^*(M+\frac{M^2}{2})
%\end{pmatrix},
%\end{align*}
being $\det(L_2) = 9u_1^*u_2^*(2+\delta_c)k_c^4\delta_c$,
%\begin{align*}
%\det(L_2) &= 9u_1^*u_2^*(2\delta_c+\delta_c^2)k_c^4,
%\end{align*}
$T_i = T^{(i)}_0+T^{(i)}_1 k_c^2 $, for $i=1,2$, and
\begin{align*}
& S_1 	=  2+\frac{Mb}{2}+2MM^*+k_c^2(\eps_M+\delta_c),\quad
 S_2 = \frac{b}{2}+2M^*\eps_M+k_c^2 M^* (\eps_M+M\delta_c),\\
& T^{(1)}_0= \frac{1}{2}(2+\frac{Mb}{2}+2MM^*), \quad T^{(1)}_1=\frac{1}{2}(1-M+2MM^*+\delta_c), \\
& T^{(2)}_0= \frac{b}{4}+M^*\eps_M,\quad T^{(2)}_1= 1+\frac{M^*}{2}(M-1+M\delta_c).
\end{align*}

\end{lemma}

\no\textbf{End of the proof of Theorem~\ref{th:example}. }
\begin{lemma}
\label{lemma:orders3}
There exists a constant $C<0$ such that $	k_c^2\delta_c\langle \bG_2,\boldeta \rangle \to C$ as $b\to0$. Consequently, 
$\langle \bG_2,\boldeta \rangle \to -\infty $ as $b\to0$.
\end{lemma}
This result together with \fer{lim:G1} implies that the solution of the Stuart-Landau equation \fer{eq:sl}, given by
\begin{align}
\label{solex:sl}
A^2(t) = \frac{\sigma}{\ell}	\frac{1}{1+\big(A_0^{-2}\frac{\sigma}{\ell}-1\big)e^{-2\sigma t}},\quad A_0=A(0),
\end{align}
satisfies $A\to0$ uniformly in $(0,\infty)$ as $b\to0$,	
which proves the result. 
$\Box$

\bigskip

\no\textbf{Proof of Lemma~\ref{lemma:orders3}. }
We set 
\begin{align*}
\zeta k_c^2\delta_c(T_1 \bv_{22}^{(1)}+ T_2 \bv_{22}^{(2)}) = I_4k_c^4 +I_2 k_c^2+I_0+I_{-2}k_c^{-2},	
\end{align*}
with $\zeta=9u_1^*u_2^*(2+\delta_c)$, $I_4  = 	\bw_2^{(1)}T_1^{(1)}+\bw_2^{(2)}T_1^{(2)}$,
$I_{-2}= \bw_0^{(1)}T_0^{(1)}+\bw_0^{(2)}T_0^{(2)}$, 
\begin{align*}
I_2 &= \bw_2^{(1)}T_0^{(1)}+\bw_1^{(1)}T_1^{(1)}+\bw_2^{(2)}T_0^{(2)}+\bw_1^{(2)}T_1^{(2)},\\
I_0 &= \bw_1^{(1)}T_0^{(1)}+\bw_0^{(1)}T_1^{(1)}+\bw_1^{(2)}T_0^{(2)} +\bw_0^{(2)}T_1^{(2)}.
\end{align*}
We have, as $b\to0$, $
\bw_0 \to (-1,3/2)$ and $ (T_0^{(1)},T_0^{(2)})\to (3/2,0)	
$, 
implying that $I_{-2}k_c^{-2}\to0$ as $b\to0$. For $I_4$, we have 
\begin{align*}
\frac{1}{4}I_4 &= 	 \Big((-u_2^*+u_1^*M)\eps_M+(-u_2^*+u_1^*M^2)\delta_c \Big) \Big(\frac{M^*-1}{2}\eps_M+\frac{1-M^*M}{2}\delta_c\Big) +O(\eps_M\delta_c)+O(\delta_c^2).
\end{align*}
Since $k_c^4\delta_c^2\to0$, $k_c^4\delta_c\eps_M\to0$ and $k_c^4\eps_M^2\to1$, we deduce
\begin{align*}
\lim_{b\to0}I_4k_c^4 = 	4\lim_{b\to0}(-u_2^*+u_1^*M)\frac{M^*-1}{2}=9.
\end{align*}
For $I_2$, we have
\begin{align*}
\bw_1^{(1)}T_1^{(1)}+\bw_1^{(2)}T_1^{(2)} &= 	\Big(-2u_2^*+4u_1^*M\Big(1+\frac{M}{2}\Big)-u_2^*bM\Big) \Big(\frac{M^*-1}{2}\eps_M+\frac{1-M^*M}{2}\delta_c\Big)\\
& +\big(-u_2^*T_1^{(1)}+(2u_2^*-u_1^*M)T_1^{(2)}\big)\eps_M +O(\delta_c)
\end{align*}
and 
$
 \bw_2^{(1)}T_0^{(1)}+\bw_2^{(2)}T_0^{(2)}	 
=(-u_2^*+u_1^*M)\eps_M (1-M^*)
+O(\delta_c)
$. 
Since $k_c^2\delta_c\to0$ and $k_c^2\eps_M\to-1$, we have
\begin{align*}
\lim_{b\to0}I_2k_c^2 = 	-\lim_{b\to0} \Big[\Big(-2u_2^*&+4u_1^*M\Big(1+\frac{M}{2}\Big)\Big) \frac{M^*-1}{2} +(-u_2^*T_1^{(1)}+(2u_2^*-u_1^*M)T_1^{(2)})\\
&+(-u_2^*+u_1^*M) (1-M^*)\Big] = \frac{5}{2}.
\end{align*}
Finally, for $I_0$, we have 
\begin{align*}
\bw_0^{(1)}T_1^{(1)}+\bw_0^{(2)}T_1^{(2)} &= 	\Big(-\frac{u_2^*}{2} +O(b)\Big) T_1^{(1)} +
\Big(u_2^*-u_1^*(M+\frac{M^2}{2}) +O(b)\Big)T_1^{(2)}
\end{align*}
\begin{align*}
\bw_1^{(1)}T_0^{(1)}+\bw_1^{(2)}T_0^{(2)} &= 	\Big(-2u_2^*+4u_1^*M\Big(1+\frac{M}{2}\Big)\Big) (T_0^{(1)}-T_0^{(2)})+O(b)+O(\eps_M) +O(\delta_c)\\
& = O(b)+O(\eps_M) +O(\delta_c).
\end{align*}
Therefore
\begin{align*}
\lim_{b\to0}I_0 = 	\lim_{b\to0} \Big(-\frac{u_2^*}{2}  T_1^{(1)} +
\Big(u_2^*-u_1^*(M+\frac{M^2}{2}) \Big)T_1^{(2)}\Big)=\frac{9}{4}
\end{align*}
implying $\lim_{b\to 0} k_c^2\delta_c(T^1 \bv_{22}^{(1)}+ T^2 \bv_{22}^{(2)}) = 55/288$.
$\Box$

%%%%%%%%%%%%%%%%%%%%%%%%%%%%%%%%%%%%%%%%%%%%%%%%%%%%%%%%%%%%%%%%%%%%%%%%%%%%%%%%%%%%

%\newpage

%\bigskip 

%\no\textbf{References}

\end{document}